\title{Basic properties of SLE}
\author{Steffen Rohde \and Oded Schramm}
\documentclass[12pt]{article}
\usepackage{amsmath}
\usepackage{amsthm}
\usepackage{amsfonts}
\usepackage{graphicx}

\newif\ifdraft
\drafttrue
\def\note#1/{\ifdraft {\bf [#1]}\fi}
\numberwithin{equation}{section}
\numberwithin{figure}{section}

\newtheorem{theorem}{Theorem}
\numberwithin{theorem}{section}
\newtheorem{corollary}[theorem]{Corollary}
\newtheorem{lemma}[theorem]{Lemma}

\newtheorem{proposition}[theorem]{Proposition}
\newtheorem{conjecture}[theorem]{Conjecture}
\newtheorem{problem}[theorem]{Problem}
\newtheorem{exercise}[theorem]{Exercise}
\theoremstyle{definition}\newtheorem*{update}{Update}
\theoremstyle{definition}\newtheorem{remark}[theorem]{Remark}

\def\eref#1{(\ref{#1})}
\def\QED{\qed\medskip}

\newcommand{\R}{\mathbb{R}}
\newcommand{\C}{\mathbb{C}}
\newcommand{\Z}{\mathbb{Z}}
\newcommand{\N}{\mathbb{N}}

\def\H{\mathbb{H}}
\def\U{\mathbb{U}}
\def\diam{\mathrm{diam}}
\def\area{\mathrm{area}}
\def\dist{\mathrm{dist}}
\def\ceil#1{\lceil{#1}\rceil}

\def\Im{{\rm Im}}
\def\Re{{\rm Re}}
\def\SLEkk#1/{$\mathrm{SLE}_{#1}$}
\def\SLEk/{\SLEkk{\kappa}/}
\def\SLE/{$\mathrm{SLE}$}
\def\Ito/{It\^o}
\def \eps {\epsilon}
\def \P {{\bf P}}
\def\md{\mid}
\def\Bb#1#2{{\def\md{\bigm| }#1\bigl[#2\bigr]}}
\def\BB#1#2{{\def\md{\Bigm| }#1\Bigl[#2\Bigr]}}
\def\Bs#1#2{{\def\md{\mid}#1[#2]}}
\def\Pb{\Bb\P}
\def\Eb{\Bb\E}
\def\PB{\BB\P}
\def\EB{\BB\E}
\def\Ps{\Bs\P}
\def\Es{\Bs\E}

\def \p {{\partial}}
\def \E {{\bf E}}

\def\closure{\overline}
\def\ev#1{{\mathcal{#1}}}
\def \proof {{ \medbreak \noindent {\bf Proof.} }}
\def\proofof#1{{ \medbreak \noindent {\bf Proof of #1.} }}

\def\xo{\xi}
\def\field{\mathcal{F}}
\def\ix{\hat x}
\def\iu{\hat u}
\def\iv{\hat v}
\def\iy{\hat y}
\def\iz{\hat z}

\def\hxo{\hat\xo}
\def\hf{\hat f}
\def\hh{\vartheta}

\def\hG{\hat G}

\def\hF{\hat F}
\def\bF{\bar F}
\def\Li{\mathrm{Li}}

\def\sign{\mathrm{sign}}
\def\BCM{MR1815718}
\def\BMakFine{MR2000g:30018}
\def\BJonMak{MR96k:30027}
\def\BDurStoch{MR97k:60148}
\def\BSmirnovPerc{MR1851632}
\def\BDup{MR2001c:82040}
\def\BDupSal{MR88d:82073}
\def\WernerSurvey{MR1905353}
\def\BhigherFunc{MR15:419i}
\def\BCardySurvey{math-ph/0103018}
\def\BLSWi{MR2002m:60159a}
\def\BLSWii{MR2002m:60159b}
\def\BLSWiii{MR1899232}
\def\BLSWlesl{math.PR/0112234}
\def\BLSWrest{MR1992830}
\def\BSchSLE{MR1776084}
\def\BPomm{MR95b:30008}
\def\BAhlfors{MR50:10211}
\def\BKR{MR98k:28004}
\def\BMR{MRlow}
\def\BRHAdd{MR98e:82018}
\def\BWeilWil{WeilWil}
\def\BBeffaraHaus{math.PR/0211322}
\def\DubedatDuality{math.PR/0303128}
\def\RevuzYor{MR2000h:60050}
\def\SSharmonicExplorer{math.PR/0310210}

\begin{document}
\maketitle

\begin {abstract}
\SLEk/ is a random growth process based on Loewner's
equation with driving parameter
a one-dimensional Brownian motion running with speed $\kappa$.
This process is intimately connected with scaling limits
of percolation clusters and with the outer boundary of Brownian motion,
and is conjectured to correspond to scaling limits of several other
discrete processes in two dimensions.

The present paper attempts a first systematic study of \SLE/.
It is proved that for all $\kappa\ne 8$ the \SLE/ trace
is a path; for $\kappa\in[0,4]$ it is a simple path; for
$\kappa\in(4,8)$ it is a self-intersecting path; and for
$\kappa>8$ it is space-filling.

It is also shown that the Hausdorff dimension of the \SLEk/
trace is a.s.\ at most $1+\kappa/8$ and that the expected
number of disks of size $\eps$ needed to cover it inside
a bounded set is at least $\eps^{-(1+\kappa/8)+o(1)}$
for $\kappa\in[0,8)$ along some sequence $\eps\searrow 0$.
Similarly, for $\kappa\ge 4$,
the Hausdorff dimension of the outer boundary 
of the \SLEk/ hull is a.s.\ at most $1+2/\kappa$, and
the expected number of disks of radius $\eps$ needed to cover
it is at least $\eps^{-(1+2/\kappa)+o(1)}$ for a sequence
$\eps\searrow0$.
\end {abstract}

\section{Introduction}

Stochastic Loewner evolution (\SLE/) is a random
process of growth of a set $K_t$.
The evolution of the set over time is
described through the normalized conformal map
$g_t=g_t(z)$ from the complement of $K_t$.  The
map $g_t$ is the solution of Loewner's differential
equation with driving parameter a one-dimensional
Brownian motion.  \SLE/, or \SLEk/,
has one parameter $\kappa\ge 0$,
which is the speed of the Brownian motion.
A more complete definition appears in Section~\ref{fund} below.

\begin{figure}
\centerline{\includegraphics*[height=2.3in]{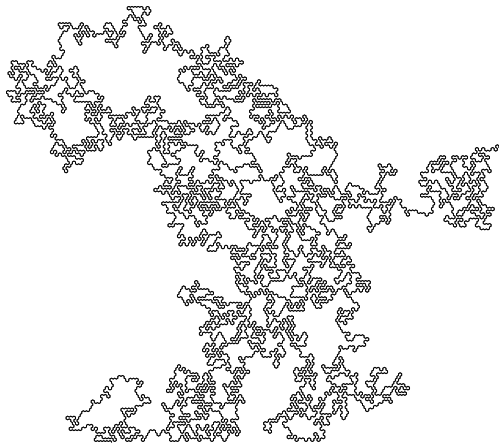}}
\caption{\label{f.percbd}The boundary of a percolation
cluster in the upper half plane, with appropriate boundary
conditions.  It converges to the chordal \SLEkk6/ trace.
}
\end{figure}

The \SLE/ process was introduced in~\cite{\BSchSLE}.
There, it was shown that under the assumption of the existence and conformal
invariance of the scaling limit of loop-erased random walk,
the scaling limit is \SLEkk2/. See Figure~\ref{f.lerw}.
It was also
stated there without proof that \SLEkk6/ is the scaling limit
of the boundaries of critical percolation clusters, assuming
their conformal invariance. Smirnov~\cite{\BSmirnovPerc} has
recently proved the conformal invariance conjecture for
critical percolation on the triangular grid
and the claim that \SLEkk6/ describes the limit.
See Figure~\ref{f.percbd}.
With the proper setup, the outer boundary of \SLEkk6/ is
the same as the outer boundary of planar Brownian
motion~\cite{\BLSWrest} (see also~\cite{\WernerSurvey}).
\SLEkk8/ has been conjectured~\cite{\BSchSLE} to be the scaling
limit of the uniform spanning tree Peano curve (see Figure~\ref{f.peano}),
and there are various further conjectures for other parameters.  Most
of these conjectures are described in Section~\ref{further} below.
Also related is the work of Carleson and Makarov~\cite{\BCM},
which studies growth processes
motivated by DLA via Loewner's equation.

\SLE/ is amenable to computations.  In~\cite{\BSchSLE} a
few properties of \SLE/ have been derived; in particular,
the winding number variance.
In the series of papers~\cite{\BLSWi,\BLSWii,\BLSWiii},
a number of other properties of \SLE/ have been studied.
The goal there was not to investigate \SLE/ for its own
sake, but rather to use \SLEkk6/ as a means for the determination
of the Brownian motion intersection exponents.

As the title suggests, the goal of the present paper is to study
the fundamental properties of \SLE/.
There are two main variants of \SLE/, {\bf chordal} and {\bf radial}.
For simplicity, we concentrate on chordal \SLE/; however, all the
main results of the paper carry over to radial \SLE/ as well.
In chordal \SLE/, the set $K_t$, $t\ge 0$, called the {\bf \SLE/ hull},
is a subset of the closed upper half plane $\closure{\H}$
and $g_t:\H\setminus K_t\to\H$ is the conformal uniformizing map,
suitably normalized at infinity.

We show that with the possible exception of $\kappa=8$, a.s.\ there
is a (unique) continuous path $\gamma:[0,\infty)\to\closure{\H}$ such that
for each $t>0$ the set $K_t$ is the union of $\gamma[0,t]$
and the bounded connected components of $\closure{\H}\setminus\gamma[0,t]$.
The path $\gamma$ is called the \SLE/ trace.  It is shown that
$\lim_{t\to\infty}|\gamma(t)|=\infty$ a.s.
We also describe two phase transitions for the \SLE/ process.
In the range $\kappa\in[0,4]$, a.s.\ $K_t=\gamma[0,t]$
for every $t\ge 0$ and $\gamma$ is a simple path.
For $\kappa\in(4,8)$ the path $\gamma$ is not a simple path
and for every $z\in\H$ a.s.\ $z\notin\gamma[0,\infty)$
but $z\in\bigcup_{t>0} K_t$.  Finally, for $\kappa>8$ we have
$\closure{\H}=\gamma[0,\infty)$ a.s.
The reader may wish to examine Figures~\ref{f.lerw},~\ref{f.percbd}
and~\ref{f.peano}, to get an idea of what the \SLEk/ trace looks
like for $\kappa=2$, $6$ and $8$, respectively.

We also discuss the expected number of disks needed to cover the \SLEk/ trace
and the outer boundary of $K_t$. It is proved that the Hausdorff dimension of
the trace is a.s.\ at most $1+\kappa/8,$ and that the Hausdorff dimension of
the outer boundary $\p K_t$ is a.s.\ at most $1+2/\kappa$ if $\kappa\ge 4$.
For $\kappa\in[0,8),$ we also show that the expected number of disks of
size $\eps$ needed to cover the trace inside a bounded set is at least
$\eps^{-(1+\kappa/8)+o(1)}$ along some sequence $\eps\searrow 0$.
Similarly, for $\kappa\ge 4,$ the expected number of disks of radius
$\eps$ needed to cover the outer boundary
is at least $\eps^{-(1+2/\kappa)+o(1)}$ for a sequence
of $\eps\searrow0$.
Richard Kenyon has earlier made the conjecture that the
Hausdorff dimension of the outer boundary is a.s.\ $1+2/\kappa$.
These results offer strong support for this conjecture.

It is interesting to compare our results to recent results for the deterministic
Loewner evolution, i.e., the solutions to the Loewner equation with a deterministic
driving function $\xo(t)$. In~\cite{\BMR} it is shown that if $\xo$ is H\"older
continuous with exponent $1/2$ and small norm, then $K_t$ is a simple path. On the
other hand, there is a function
$\xo,$ H\"older continuous with exponent $1/2$ and having large norm, such that
$K_t$ is not even locally connected,
and therefore
there is no continuous path $\gamma$ generating $K_t$.
In this example, $K_t$ spirals infinitely often around a disk $D$, accumulating on
$\partial D$, and then spirals out again.  It is easy to see that the disk $D$ can be
replaced by any compact connected subset of $\H$.
Notice that according to the law of the iterated logarithm,
a.s.\ Brownian motion is not H\"older continuous with exponent $1/2$. Therefore, it seems
unlikely that the results of the present paper can be obtained from deterministic results.

Our results are based on the computation and estimates of
the distribution of $|g_t'(z)|$ where $z\in\H$.
Note that in \cite{\BLSWii} the derivatives $g_t'(x)$ are studied for $x\in\R$.

The organization of the paper is as follows.  Section~\ref{fund} introduces the 
basic definitions and some fundamental properties.  The goal of Section~\ref{Sderexp}
is to obtain estimates for quantities related to $\Eb{|g_t'(z)|^a}$, for various 
constants $a$ (another result of this nature
is Lemma~\ref{howclose}),
and to derive some resulting continuity properties of $g_t^{-1}$.
Section~\ref{Reduction} proves a general criterion for the existence of
a continuous trace, which does not involve randomness.  The
proof that the \SLEk/ trace is continuous for $\kappa\ne8$ is
then completed in Section~\ref{scont}.  There, it is also proved
that $g_t^{-1}$ is a.s.\ H\"older continuous when $\kappa\ne 4$.
Section~\ref{sphases} discusses the two phase transitions
$\kappa=4$ and $\kappa=8$ for \SLEk/.  Besides some quantitative properties,
it is shown there that the trace is a.s.\ a simple path iff $\kappa\in[0,4]$,
and that the trace is space-filling for $\kappa>8$.
The trace is proved to be transient when $\kappa\ne8$ in Section~\ref{trans}.
Estimates for the dimensions of the trace and the boundary of the hull
are established in Section~\ref{dimensions}.  Finally, a collection of open problems
is presented in Section~\ref{further}.

\begin{update}
Since the completion and distribution of the first version of this paper,
there has been some further progress.  In \cite{\BLSWlesl} it was proven that
the scaling limit of loop-erased random walk is \SLEkk2/ and the
scaling limit of the UST Peano path is \SLEkk8/.  As a corollary of
the convergence of the UST Peano path to \SLEkk8/, it was also established
there that \SLEkk8/ is generated by a continuous transient path, thus 
answering some of the issues left open in the current paper.  However,
it is quite natural to ask for a more direct analytic proof of these 
properties of \SLEkk8/.

Recently, Vincent Beffara~\cite{\BBeffaraHaus} has announced a proof
that the Hausdorff dimension of the \SLEk/ trace is
$1+\kappa/8$ when $4\ne\kappa\le 8$.

The paper~\cite{\SSharmonicExplorer} proves the convergence of the
harmonic explorer to \SLEkk4/.
\end{update}

\section{Definitions and background}\label{fund}
\subsection{Chordal \SLE/}

Let $B_t$ be Brownian motion on $\R$, started from $B_0=0$.
For $\kappa\ge 0$ let $\xo(t):=\sqrt{\kappa} B_t$ and
for each $z\in\closure{\H}\setminus\{0\}$ let
$g_t(z)$ be the solution of the ordinary differential equation
\begin{equation}\label{chordal}
\p_t g_t(z)= \frac 2{g_t(z)-\xo(t)}\,,\qquad g_0(z)=z.
\end{equation}
The solution exists as long as $g_t(z)-\xo(t)$ is bounded
away from zero.  We denote by $\tau(z)$ the first time
$\tau$ such that $0$ is a limit point of $g_t(z)-\xo(t)$
as $t\nearrow\tau$.  Set
$$
H_t:= \bigl\{z\in\H:\tau(z)>t\bigr\}\,,
\qquad
K_t:= \bigl\{z\in\closure{\H}:\tau(z)\le t\bigr\}\,.
$$
It is immediate to verify that $K_t$ is compact and
$H_t$ is open for all $t$.
The parameterized collection of maps $(g_t:t\ge 0)$
is called {\bf chordal \SLEk/}.
The sets $K_t$ are the {\bf hulls} of the \SLE/.
It is easy to verify that for every $t\ge 0$ the map
$g_t:H_t\to\H$ is a conformal
homeomorphism and that $H_t$ is the unbounded component
of $\H\setminus K_t$.
The inverse of $g_t$ is obtained by flowing backwards
from any point $w\in \H$ according to the equation~\eref{chordal}.
(That is, the fact that $g_t$ is invertible is a particular case of
a general result on solutions of ODE's.)
One only needs to note that in this backward flow,
the imaginary part increases, hence the point cannot hit the singularity.
It also cannot escape to infinity in finite time.
The fact that $g_t(z)$ is analytic in $z$ is clear,
since the right hand side of~\eref{chordal} is analytic in $g_t(z)$.

The map $g_t$ satisfies the so-called {\bf hydrodynamic} normalization
at infinity:
\begin{equation}\label{hydro}
\lim_{z\to\infty} g_t(z)-z=0\,.
\end{equation}
Note that this uniquely determines $g_t$ among conformal
maps from $H_t$ onto $\H$.
In fact,~\eref{chordal} implies that $g_t$ has the
power series expansion
\begin{equation}
\label{power}
g_t(z)=z+ \frac{2t}z +\cdots,\qquad z\to\infty\,.
\end{equation}
The fact that $g_t\ne g_{t'}$ when $t'>t$  implies that
$K_{t'}\ne K_t$.  The relation $K_{t'}\supset K_t$ is clear.

Two important properties of chordal \SLE/ are scale-invariance
and a sort of stationarity.  These are summarized in the following
proposition. (A similar statement appeared in~\cite{\BLSWi}.)

\begin{proposition}\label{funda}
\par\noindent{\rm(i)}
\SLEk/ is scale-invariant in the following
sense.  Let $K_t$  be the hull of \SLEk/, and let $\alpha>0$.
Then the process $t\mapsto\alpha^{-1/2} K_{\alpha t}$
has the same law as $t\mapsto K_t$.  The process
$(t,z)\mapsto\alpha^{-1/2}\, g_{\alpha t}(\sqrt{\alpha}z)$ has the same law
as the process $(t,z)\mapsto g_t (z)$.
\par\noindent{\rm(ii)}
Let $t_0>0$.  Then the map
$(t,z)\mapsto \tilde g_t(z):=
g_{t+t_0} \circ g_{t_0}^{-1}\bigl(z+\xo(t_0)\bigr)-\xo(t_0)$
has the same law as $(t,z)\mapsto g_t(z)$;
moreover, $(\tilde g_t)_{t\ge 0}$ is independent from $(g_{t})_{0\le t\le t_0}$.
\end{proposition}

The scaling property easily follows from the scaling property
of Brownian motion, and the second property follows from the
Markov property and translation invariance of Brownian motion.
One just needs to write down the expression for $\p_t \tilde g_t$.
We leave the details as an exercise to the reader.

\medskip
The following notations will be used throughout the paper.
$$
f_t:= g_t^{-1},\qquad
\hf_t(z):= f_t\bigl(z+\xo(t)\bigr).
$$
The {\bf trace} $\gamma$ of \SLE/ is defined by
$$
\gamma(t):= \lim_{{z\to 0}} \hf_t(z)\,,
$$
where  $z$ tends to $0$ within $\H$.  If the limit
does not exist, let $\gamma(t)$ denote the set of all limit points.
We say that the \SLE/ trace is a continuous path
if the limit exists for every $t$ and $\gamma(t)$ is
a continuous function of $t$.

\subsection{Radial \SLE/}

Another version of \SLEk/ is called {\bf radial \SLEk/}.
It is similar to chordal \SLE/ but appropriate for the situation where there
is a distinguished point in the interior of the domain.
Radial \SLEk/ is defined as follows.  Let $B(t)$ be Brownian motion
on the unit circle $\p\U$, started from a uniform-random point
$B(0)$, and set $\xo(t):= B(\kappa t)$.  The conformal maps $g_t$ are defined
as the solution of
\begin{equation*}
\p_t g_t(z)= -g_t(z)\,\frac {g_t(z)+\xo(t)}{g_t(z)-\xo(t)}\,,\qquad g_0(z)=z\,,
\end{equation*}
for $z\in\U$.
The sets $K_t$ and $H_t$ are defined as for chordal \SLE/.
Note that the scaling property~\ref{funda}.(i) fails for radial \SLE/.
Mainly due to this reason, chordal \SLE/ is easier to work with.
However, appropriately stated, all the main results of this paper are valid
also for radial \SLE/.
This follows from~\cite[Prop.~4.2]{\BLSWii}, which says in a precise
way that chordal and radial \SLE/ are equivalent.

\subsection{Local martingales and martingales}

The purpose of this subsection is to present a slightly technical lemma
giving a sufficient condition for a local martingale to be a martingale.
Although we have not been able to find an appropriate reference, the
lemma must be known (and is rather obvious to the experts).

See, for example, \cite[\S IV.1]{\RevuzYor} for a discussion
of the distinction between
a local martingale and a martingale.

While the stochastic calculus needed for the rest of the paper is
not much more than familiarity with \Ito/'s formula, this subsection
does assume a bit more.

\begin{lemma}\label{l.locm}
Let $B_t$ be stardard one dimensional Brownian motion, and let
$a_t$ be a progressive real valued locally bounded process.  Suppose that
$X_t$ satisfies
$$
X_t = \int_0^t a_s\,dB_s\,,
$$
and that for every $t>0$ there is a finite constant $c(t)$ such that
\begin{equation}\label{e.as}
a_s^2 \le c(t) X_s^2+c(t)
\end{equation}
for all $s\in[0,t]$ a.s.  Then $X$ is a martingale.
\end{lemma}

\proof
We know that $X$ is a local martingale.
Let $M>0$ be large, and let $T:=\inf\{t:|X_t|\ge M\}$.
Then $Y_t:=X_{t\wedge T}$ is a martingale (where $t\wedge T=\min\{t,T\}$).
Let $f(t):=\Eb{Y_t^2}$.  
\Ito/'s formula then gives
$$
f(t')=\EB{ \int_0^{t'} a_s^2\,1_{s<T}\,ds}\,.
$$
Our assumption $a_s^2\le c(t)X_s^2+c(t)$ therefore implies
that for $t'\in[0,t]$
\begin{equation}\label{e.intiq}
f(t')\le c(t)\,t'+ c(t) \int_0^{t'} f(s)\,ds\,.
\end{equation}
This implies $f(s)< \exp(2\, c(t)\, s)$ for all $s\in[0,t]$,
since~\eref{e.intiq} shows that $t'$ cannot
be the least $s\in[0,t]$ where $f(s)\ge \exp(2\, c(t)\, s)$.
Thus, 
$$
\Eb{\left<X,X\right>_{t\wedge T}}=
\Eb{\left<Y,Y\right>_t}
=\Eb{Y_t^2}=f(t)<\exp(2\, c(t)\, t)\,.
$$
Taking $M\to\infty$, we get by
monotone convergence $\Eb{\left<X,X\right>_t}\le \exp(2\, c(t)\, t)<\infty$.
Thus $X$ is a martingale (for example, by~\cite[IV.1.25]{\RevuzYor}).
\QED

\section{Derivative expectation}\label{Sderexp}

In this section, $g_t$ is the \SLEk/ flow; that is,
the solution of~\eref{chordal} where
$\xo(t):=B(\kappa t)$, and $B$ is standard Brownian motion on $\R$
starting from $0$.  Our only assumption on $\kappa$ is
$\kappa>0$.  The goal of the section is to derive bounds
on quantities related to $\Eb{|g_t'(z)|^a}$.  Another
result of this nature is Lemma~\ref{howclose},
which is deferred to a later section.

\subsection{Basic derivative expectation}

We will need estimates for the moments of $|\hf_t'|$. In this subsection, we will
describe a change of time and obtain derivative estimates for the changed time.

For convenience, we take $B$ to be two-sided Brownian motion.
The equation~\eref{chordal} can also be solved for negative
$t$, and $g_{t}$ is a conformal map from $\H$ into a subset
of $\H$ when $t<0$. Notice that the scale invariance (Proposition~\ref{funda}.(i))
also holds for $t<0$.

\begin{lemma}\label{negt}
For all fixed $t\in\R$ the map
$z\mapsto g_{-t}\bigl(z\bigr)$ has the
same distribution as the map $z\mapsto \hf_t(z)-\xo(t)$.
\end{lemma}
\proof
Fix $t_1\in\R$, and let
\begin{equation}\label{hxo}
\hxo_{t_1}(t)=\xo(t_1+t)-\xo(t_1)\,.
\end{equation}
Then $\hxo_{t_1}:\R\to\R$ has the same law as $\xo:\R\to\R$.
Let
$$
\hat g_t(z):=g_{t_1+t}\circ g_{t_1}^{-1}\bigl(z+\xo(t_1)\bigr)-\xo(t_1)\,,
$$
and note that $\hat g_0(z)=z$ and
$\hat g_{-t_1}(z)=\hf_{t_1}(z)-\xo(t_1)$.
Since
$$
\p_t\hat g_t = \frac 2{\hat g_t+\xo(t_1)-\xo(t+t_1)}
=
 \frac 2{\hat g_t-\hxo_{t_1}(t)}\,,
$$
the lemma follows from~\eref{chordal}.
\QED

Note that~\eref{chordal} implies that $\Im\bigl(g_t(z)\bigr)$
is monotone decreasing in $t$ for every $z\in\H$.
For $z\in\H$ and $u\in\R$ set
\begin{equation}\label{Tdef}
T_u=T_u(z):=\sup\Bigl\{t\in\R: \Im\bigl(g_{t}(z)\bigr)\ge e^u\Bigr\}\,. %
\end{equation}
We claim that for all $z\in\H$ a.s.\ $T_u\ne\pm\infty$.  
By~\eref{chordal}, $|\p_t g_t(z)| = 2 |g_t(z)-\xo(t)|^{-1}$.
Setting $\bar\xo(t):=\sup\bigl\{|\xo(s)|:s\in[0,t]\bigr\}$,
the above implies that 
$ |g_t(z)|\le |z|+ \bar\xo(t) + 2 \sqrt{t} $ for $t<\tau(z)$,
since $\p_t|g_t(z)|\le 2/(|g_t(z)|-\bar\xo(t))$ whenever
$|g_t(z)|>\bar\xo(t)$. Setting $y_t:=\Im g_t(z)$, we get from~\eref{chordal},
$$
-\p_t \log y_t \ge 2\,(|z|+2\,\bar\xo(t)+2\sqrt t)^{-2}\,.
$$
The law of iterated logarithms implies that the right hand side is
not integrable over $[0,\infty)$ nor over $(-\infty,0]$.  Thus,
$|T_u|<\infty$ a.s.
\medskip

We will need the formula
\begin{equation}\label{theder}
\begin{aligned}
\p_t \log|g_t'(z)|
&= \Re\Bigl(\frac{\p_z\p_t g_t(z)}{g_t'(z)}\Bigr)
=\Re\Bigl(g_t'(z)^{-1}\,\p_z\frac{2}{g_t(z)-\xo(t)}\Bigr)
\\ &
=-2\,\Re\Bigl( \bigl(g_t(z)-\xo(t)\bigr)^{-2}\Bigr)
\,.
\end{aligned}
\end{equation}

Set $u=u(z,t):=\log\Im g_t(z)$.
Observe that~\eref{chordal} gives
\begin{equation}\label{du}
\p_t u= -2\, |g_t(z)-\xo(t)|^{-2}\,,
\end{equation}
and~\eref{theder} gives
\begin{equation}\label{dlog}
\p_u \log \bigl| g_t'(z)\bigr|
=\frac{\Re\bigl((g_t(z)-\xo(t))^2\bigr)}{\left|g_t(z)-\xo(t)\right|^2}\,.
\end{equation}
\medskip

Fix some $\iz=\ix +i\iy\in\H$.  For every $u\in\R$, let
$$
z(u):= g_{T_u(\iz)}(\iz)-\xo(T_u),\qquad
\psi(u):= \frac{\iy}{ y(u)}\, \bigl|g_{T_u}'(\iz)\bigr|,
$$
and
$$
x(u):= \Re(z(u)),\qquad
y(u):=\Im(z(u))=\exp(u).
$$

\begin{theorem}\label{derexp}
Let $\iz=\ix+i\iy\in\H$ as above.
Assume that $\iy\ne 1$,
and set $\nu:=-\mathrm{sign}(\log \iy)$. 
Let $b\in\R$.  Define $a$ and $\lambda$ by
\begin{equation}\label{al}
a:=2\, b+ \nu\,\kappa\, b\,(1-b)/2,
\qquad
\lambda:= 4\, b+ \nu\,\kappa\,b \,(1-2 b)/2\,.
\end{equation}
Set
$$
F(\iz)=F_{b}(\iz) :=
\iy^a\,\EB{\bigl(1+x(0)^2\bigr)^b\,\bigl|g_{T_0(\iz)}'(\iz)\bigr|^a}.
$$
Then
\begin{equation}\label{Feq}
F(\iz) =
\Bigl(1+(\ix/\iy)^2\Bigr)^b\iy^{\lambda}\,.
\end{equation}
\end{theorem}

Before we give the short proof of the theorem, a few remarks may be of
help to motivate the formulation and the proof. 
Our goal was to find an expression for
$\EB{\bigl|g_{T_0(\iz)}'(\iz)\bigr|^a}$. It turns out to be more convenient
to consider
$$
\bF(\iz):=\iy^a\,\EB{\bigl|g_{T_0(\iz)}'(\iz)\bigr|^a}
= \EB{\psi(0)^a}.
$$
The obvious strategy is to find a differential equation which $\bF$ must
satisfy and search for a solution.  The first part is not too difficult,
and proceeds as follows.
Let $\field_u$ denote the $\sigma$-field generated
by $\left<\xo(t): (t-T_u)\nu\ge 0\right>$.
Note that the strong Markov property for $\xo$ and the chain rule imply
that for $u$ between $0$ and $\iu:=\log \iy$,
\begin{equation}\label{mart}
\EB{\psi(0)^a\md\field_u}= \psi(u)^a\,\bF\bigl(z(u)\bigr)\,.
\end{equation}
Hence, the right hand side is a martingale.
Observe that
\begin{equation}\label{chng}
dx = \frac{2\, x\,dt}{ x^2+y^2}-d\xo,
\qquad
dy = \frac{ -2\, y\,dt}{ x^2+y^2},
\qquad
d\log\psi = \frac{4\, y^2\,dt}{ (x^2+y^2)^2}\,.
\end{equation}
(The latter easily follows from~\eref{theder} and~\eref{du}.)
We assume for now that $\bF$ is smooth.
\Ito/'s formula may then be calculated for the right hand side
of~\eref{mart}. Since it is a martingale,
the drift term of the It\^o derivative must vanish; that is,
$$
\psi^a\cdot\Lambda \bF=0\,,
$$
where
\begin{equation}\label{dLam}
\Lambda G:=
 -\frac{4\,\nu\,a\, y^2}{ (x^2+y^2)^2} \,G
-\frac{2\,\nu\, x}{ x^2+y^2}\, \p_x G
+\frac{2\,\nu\, y}{ x^2+y^2}\, \p_y G
+\frac{\kappa}{ 2}\, \p_x^2 G\,.
\end{equation}
(The $-\nu$ factor comes from the fact that
$t$ is monotone decreasing with respect to the
filtration $\field_u$ iff $\nu=1$.) 
Guessing a solution for the equation $\Lambda G=0$ is not too
difficult (after changing to coordinates where scale invariance is
more apparent).  It is easy to verify that
$$
\hF(x+iy)=\hF_{b,\,\lambda}(x+iy):=
\Bigl(1+(x/y)^2\Bigr)^by^{\lambda}\,,
$$
satisfies $\Lambda\hF =0$.
Unfortunately, $\hF$ does not satisfy the boundary values
$\hF=1$ for $y=1$, which hold for $\bF$.
Consequently, the theorem gives a formula for $F$, rather than
for $\bF$.
(Remark~\ref{r.imp} concerns the problem of determining $\bF$.)
  Assuming that $F$ is $C^2$,
the above derivation does apply to $F$, and shows
that $\Lambda F=0$.  However, we have not found a clean reference to the
fact that $F\in C^2$.
Fortunately, the proof below does not need to assume this.

\proofof{Theorem \ref{derexp}}
Note that by~\eref{du}
$$
du=-2\,|z|^{-2}\,dt.
$$
Let
$$
\hat B(u):=-\sqrt {2/\kappa}\, \int_{t=0}^{T_u}|z|^{-1}\, d\xo.
$$
Then $\hat B$ is a local martingale and
$$
d\langle \hat B\rangle =  - 2\,|z|^{-2}\, dt = du.
$$
Consequently,
$\hat B(u)$ is Brownian motion (with respect to $u$). 
Set $M_u:=\psi(u)^a\,\hF(z(u))$.
\Ito/'s formula gives
$$
dM_u = -2 M\,\frac { b\,x}{ x^2+y^2}\,d\xo
=
 \sqrt{2\,\kappa}\, M\,\frac { b\,x}{{\sqrt{ x^2+y^2}}}\,d\hat B
\,.
$$
Thus $M$ is a local martingale.  In fact,
Lemma~\ref{l.locm} then tells us that $M$ is a martingale.
Consequently, we have 
\begin{equation}\label{martconc}
\psi(\iu)^a\,\hF(\iz)=\Eb{\psi(0)^a\, \hF(z(0))}
\end{equation}
and the theorem clearly follows.
\QED

%
%
%

In Section~\ref{dimensions} we will estimate the expected number of disks
needed to cover the boundary of $K_t$.
To do this, we need the following lower bound on
the expectation of the derivative.

\begin{lemma}\label{lbdexp}
Let $\kappa>0$ and $b<(\kappa+4)/(4\kappa)$,
and define $a$ and $\lambda$ by~\eref{al} with $\nu=1$.
Then there is a constant $c=c(\kappa,b)>0$ such that
$$
\EB{|\hf_1'(\iz)|^a\,1_{\Im(\hf_1(\iz))\ge c}}\ge c %
\Bigl(1+(\ix/\iy)^2\Bigr)^b\,\iy^{\lambda-a}
$$
holds for every $\iz=\ix+i\iy\in\H$ satisfying $|\iz|\le c$. %
\end{lemma}

\proof
As before, let $u:=\log y$.
Set $v=v(u):=\sinh^{-1}(x/y)$; that is, $x=y\sinh(v)$.
Then \Ito/'s formula gives
$$
dv=
-\frac {d\xo}{|z|} + \frac 12 \,(8+\kappa\,\nu)\,|z|^{-3}\,x\, dt\,.
$$
The formula for $dv$ in terms of $\hat B$
and $du$  is
\begin{equation}\label{dv}
dv=
\sqrt{\kappa/2}\,d\hat B - \frac x{4\,|z|}\, (8+\kappa\,\nu)\, du\,.
\end{equation}

Recall that $\psi(u)^a\,F\bigl(z(u)\bigr)$ is a martingale, since
$F=\hF$.
Define a new probability measure by
$$
\tilde \P[\ev A]:=
\frac {
\EB{1_{\ev A}\,\bigl(1+x(0)^2\bigr)^b\psi(0)^a}
}{
\EB{\bigl(1+x(0)^2\bigr)^b\psi(0)^a}
}\,,
$$
for every event $\ev A$.
This is the so-called Doob-transform (or h-transform)
corresponding to the martingale $\psi(u)^a\,F\bigl(z(u)\bigr)$.
Recall that if $\alpha(w)$ is a positive martingale for a diffusion process
$dw=q_1(w,t)\,dB(t)+q_2(w,t)\,dt$, $t\le t_1$, where $B(t)$ is
Brownian motion, then for $t<t_1$,
$dw = %
q_1(w,t)\,d\tilde B(t)+ q_2(w,t)\,dt+ q_1(w,t)\,\p_w \log \alpha(w)\,
q_1(w,t)\,dt $,
where $\tilde B$ is Brownian motion with respect to the probability
measure weighted by $\alpha(w(t_1))$; that is, the
Doob transform of $\alpha$.  This follows, for example, from
Girsanov's Formula~\cite[\S2.12]{\BDurStoch}. %
We apply this with $w=(v,\psi)$ and $u$ as the time parameter
(in this case, $q_1$ and $q_2$ are vectors and $\p_w\log\alpha(w)$
is a linear functional),
and get by~\eref{dv}
$$
dv =  \sqrt{\kappa/2}\,
d\tilde B - \frac x{4\,|z|}\, (8+\kappa)\, du+
\frac12\,\kappa\, \p_v\log F\,du\,,
$$
where $\tilde B(u)$ is Brownian motion under $\tilde \P$.
This simplifies to
\begin{equation}\label{ndv}
dv = \sqrt{\kappa/2}\,
d\tilde B -\tilde b \,\tanh (v)\,du\,,
\end{equation}
where
$$\tilde b:= 2-\bigl(b-\frac 14\bigr)\kappa>1\,.$$
Thus, $v(u)$ is a very simple diffusion process.
When $|v|$ is large, $\tanh(v)$ is close to $\sign(v)$,
and $v$ has almost constant drift $-\tilde b\,\sign(v)$,
pushing it towards $0$.

At this point, do not assume $|\iz|\le c$, but only
$|\iz|\le 1$.
Let $\Psi:[\iu,0]\to\R$ be the continuous function that is
equal to $1$ at $0$, equal to $|\iu|+1$ at
$\iu$, has slope $-\bigl(\tilde b\wedge 2+1\bigr)/2$ in the interval
$[\iu,\iu/2]$ and has constant slope in the interval
$[\iu/2,0]$, and let
$\ev A$ be the event
$$
\ev A
:=\Bigl\{
\forall u\in[\iu,0],
\quad |v(u)|\le \Psi(u)
\Bigr\}\,.
$$
Note that  our assumption $|\iz|\le 1$ implies
that $|\iv|\le  |\iu|+\log 2$.
Since $\tilde b>1$ it easily follows from~\eref{ndv}
that there is a constant
$c_1=c_1(b,\kappa)>0$ such that $\tilde \P[\ev A]>c_1$.
In particular, $c_1$ does not depend on $\iz$.
This means
$$
\EB{1_{\ev A}\, \bigl(1+x(0)^2\bigr)^b\,\psi(0)^a}\ge
c_1\,\EB{ \bigl(1+x(0)^2\bigr)^b\,\psi(0)^a}
= c_1\,
\Bigl(1+(\ix/\iy)^2\Bigr)^b\,\iy^{\lambda}\,.
$$
However, note that $v(0)$ and hence $x(0)$ are bounded
on $\ev A$.  Therefore, there is some constant $c_2>0$
such that
\begin{equation}\label{aml}
\EB{1_{\ev A}\, \bigl|g'_{T_0}(\iz)\bigr|^a}\ge
 c_2\,
\Bigl(1+(\ix/\iy)^2\Bigr)^b\,\iy^{\lambda-a}\,.
\end{equation}

We now estimate $T_0$ on the event $\ev A$.
Recall that $T_0\le 0$.
{}From~\eref{du}
we have
$ \p_u T_u= - |z|^2/2$ and therefore on $\ev A$
\begin{align*}
T_0 & = -\int_{\iu}^0 \frac {y^2}{2}\,du-
 \int_{\iu}^0 \frac {x^2}{2}\,du
=
\frac{ \iy^2}4-\frac14 -
 \frac 12 \,\int_{\iu}^0 \sinh(v)^2 \,e^{2\, u}\,du
\\&
\ge
-\frac14 -\frac12 \int_{\iu}^0 e^{2\Psi(u)+2u}\,du
\ge -c_3\,
,
\end{align*}
where $c_3=c_3(b,\kappa)<\infty$ is some constant.
That is, we have $T_0\in[-c_3,0]$ on $\ev A$.
On the event $T_0\ge -c_3$, we clearly have
$\Im\bigl(g_{-c_3}(\iz)\bigr)\ge 1$ and also
$\bigl|g_{T_0}'(\iz)\bigr|\le\bigl|g_{-c_3}'(\iz)\bigr|$,
by~\eref{theder}.
Consequently, Lemma~\ref{negt} and~\eref{aml} give
\begin{equation}\label{cc}
\EB{\bigl|\hf'_{c_3}(\iz)\bigr|^a\,1_{\Im(\hf_{c_3}(\iz))\ge 1}}
\ge %
 c_4\,
\Bigl(1+(\ix/\iy)^2\Bigr)^b\,\iy^{\lambda-a}\,.
\end{equation}

This is almost the result that we need. However, %
we want to replace $c_3$ by $1$.
For this, we apply scale invariance.  In this procedure,
the assumption $|\iz|\le 1$ needs to be replaced by
$|\iz|\le 1/\sqrt{c_3}$.  The proof of the lemma is now
complete.
\QED

It is not too hard to see that for every constant $C>0$ %
the statement of the lemma can be strengthened to allow
$c<|\iz|<C$.  The constant $c$ must then also depend on $C$.

\begin{remark}\label{r.imp} %
Suppose that we take
$$ b=\frac14 + \frac{2\,\nu}{ \kappa}\,,$$
and define $a$ and $\lambda$ using~\eref{al}.
Define $\tilde\P$ as in the proof of the lemma.
Then as the proof shows, $v$ becomes Brownian
motion times $\sqrt{\kappa/2}$ under $\tilde\P$,
since $\tilde b$ vanishes.  This allows a
detailed understanding of the distribution of
$|g'_{T_0}(\iz)|$.  For example, one can determine in this way the
decay of $\Pb{\psi(0)>1/2}$ as $\iy\searrow 0$.
It also follows that for such $a,b,\lambda$, and every $A\subset\R$
one can write down an explicit expression for
$\EB{|g_{T_0}'(\iz)|^a\,1_{x(0)\in A}}$,
since for every $v_0\in\R$ we have
$$
\Lambda\Bigl(\bigl(1+(x/y)^2\bigr)^b \, y^\lambda\, \frac 1{\sqrt{ \nu\, u}}\,
\exp\frac{(v_0-v)^2}{\nu\,\kappa\, u}\Bigr)=0\,.
$$
This equation is the PDE facet of the fact that $v$ is Brownian motion
times $\sqrt{\kappa/2}$ under $\tilde\P$.
However, these results will not be needed in the present paper.
\end{remark}

\subsection{Derivative upper bounds at a fixed time $t$}

{}From Theorem~\ref{derexp} it is not hard to obtain
estimates for $|\hf_t'|:$

\begin{corollary}\label{derest}
Let $b\in[0,1+4\,\kappa^{-1}]$, and
define $\lambda$ and $a$ by~\eref{al} with $\nu=1$.
There is a constant $C(\kappa,b)$, depending only
on $\kappa$ and $b$, such that the following estimate holds for all
$t\in[0,1]$, $y,\delta\in(0,1]$ and $x\in\R$.
\begin{equation}\label{kapcas}
\Pb{ |\hf_t'(x+i\,y)|\ge \delta\,y^{-1}}\le
C(\kappa,b)\, (1+x^2/y^2)^b\,(y/\delta)^\lambda\,\hh(\delta,a-\lambda)\,,
\end{equation}
where
\begin{equation*}
\hh(\delta,s)=\begin{cases}\delta^{-s}& s>0,\\
1+|\log \delta|& s=0,\\
1&s<0.\\
\end{cases}
\end{equation*}
\end{corollary}

\proof
Note that $a\ge 0$.
We assume that $\delta>y$, for otherwise the right hand side
is at least $C(\kappa,b)$, and we take $C(\kappa,b)\ge 1$.
Take $z=x+i\,y$.
By Lemma~\ref{negt}, $\hf_t'(z)$ has the
same distribution as $g_{-t}'(z)$.
Let $u_1:=\log\Im\bigl(g_{-t}(x+iy)\bigr)$.  Recall the
notation $T_u$ from~\eref{Tdef}.
Note that
\begin{equation}\label{nearhight}
\bigl|g_{-t}'(z)/g'_{T_u(z)}(z)\bigr|\le \exp\bigl(|u-u_1|\bigr),
\end{equation}
since $\bigl|\partial_u\log |g'_t(z)|\bigr|\le 1$, by~\eref{dlog}.
Moreover, it is clear that there is a constant $c$
such that $u_1\le c$, because $t,y\le 1$.
Consequently,
$$
\Pb{
|g_{-t}'(z)|\ge O(1) \delta\,y^{-1}}
\le O(1) \sum_{j=\lceil \log y \rceil}^0
\Pb{|g'_{T_j(z)}(z)|\ge \delta\,y^{-1}}.
$$
The Schwarz Lemma implies that
$y|g'(z)|\le \Im\bigl(g(z)\bigr)$ if $g:\H\to\H$.
Therefore, the above gives
\begin{equation}\label{ssum}
\Pb{
|g_{-t}'(z)|\ge O(1)\,\delta\,y^{-1}}
\le O(1)
\sum_{j=\lceil \log\delta \rceil}^0
\Pb{|g'_{T_j(z)}(z)|\ge \delta\,y^{-1}}.
\end{equation}
By scale invariance, we have for all $j\in[\log y, 0]$
$$
\Eb{y^a\,e^{-ja}\,|g'_{T_j(z)}(z)|^a}\le F_b\bigl(e^{-j} z\bigr)\,,
$$
where $F=F_b$ is as in Theorem~\ref{derexp}.
Hence
\begin{align*}
\Pb{|g'_{T_j(z)}(z)|\ge \delta\,y^{-1}}
&
=
\Pb{|g'_{T_j(z)}(z)|^a\, y^a\,\delta^{-a}\ge 1}
\\&
\le
\Eb{|g'_{T_j(z)}(z)|^a\, y^a\,\delta^{-a}}
\le \delta^{-a}\,e^{j\,a}\, F_b\bigl(e^{-j} z\bigr)
\,.
\end{align*}
Consequently, by~\eref{ssum} and Theorem~\ref{derexp},
\begin{equation*}
\Pb{
|g_{-t}'(z)|\ge O(1)\,\delta\,y^{-1}}
\le O(1)\,(1+x^2/y^2)^b\,
\delta^{-a}\,y^\lambda
\sum_{j=\lceil\log\delta\rceil}^0
e^{j(a-\lambda)}.
\end{equation*}
The corollary follows.
\QED

\subsection{Continuity of $\hf_t(0)$}

In the deterministic example of non-locally connected hulls described in the introduction, there is a
time $t_0$ for which $\lim_{z\to 0} \hf_{t_0}(z)$ does not exist
(the limit set is the outer boundary of the prescribed compact set).
Even when the \SLE/ trace is a continuous path,  it is not always true
that $(z,t)\mapsto f_t(z)$ extends continuously to $\overline{\H}\times[0,\infty)$
(this is only true for simple paths).
The next theorem shows that $\hf_t(0)=f_t(\xo(t))$ exists as a radial limit
and is continuous. Together with the result of Section~\ref{Reduction} below,
this is enough to show that the \SLE/ trace is a path.

\begin{theorem}\label{extend}
Define
$$
H(y,t):= \hf_t(i\,y),\qquad y>0,\ t\in[0,\infty).
$$
If $\kappa\ne 8$, then a.s.\ $H(y,t)$ extends continuously
to $[0,\infty)\times [0,\infty)$.
\end{theorem}

\proof
Fix $\kappa\ne 8$.
By scale invariance,
it is enough to show continuity of $H$ on $[0,\infty)\times [0,1)$.
Given $j,k\in\N$, with $k< 2^{2j}$, let $R(j,k)$ be the rectangle
$$
R(j,k):=[2^{-j-1},2^{-j}]\times
[k\, 2^{-2j},(k+1)\,2^{-2j}]\,,
$$
and set
$$
d(j,k):=\diam H\bigl(R(j,k)\bigr)\,.
$$
Take %
$b= (8+\kappa)/(4\kappa)$
and let $a$ and $\lambda$ be given by~\eref{al} with $\nu=1$.
Note that $\lambda>2$.
Set $\sigma_0:=(\lambda-2)/\max\{a,\lambda\}$,
and let $\sigma\in(0,\sigma_0)$.
Our objective is to prove
\begin{equation}\label{Asm}
\sum_{j=0}^\infty\sum_{k=0}^{2^{2j}-1}
\Pb{d(j,k)\ge 2^{-j\sigma} }<\infty\,.
\end{equation}

Fix some such pair $(j,k)$.
The idea of the proof is to decompose the time interval
$[k\, 2^{-2j}, (k+1)\,2^{-2j}]$ into (random) sub-intervals
using partition points
$\hat t_N, \hat t_{N-1},\dots,\hat t_0$ such that the change of
$\xo$ on each interval is controlled. Then it is not hard to
control $\hf_t$ on each sub-interval. Now, the details.
Set $t_0=(k+1)\,2^{-2j}$.
Inductively, let
$$
t_{n+1}:=\sup\bigl\{t<t_n : |\xo(t)-\xo(t_n)|= 2^{-j}\bigr\}\,,
$$
and let $N$ be the least $n\in\N$ such
that $t_N \le t_0-2^{-2j}$.
Also set $t_\infty:=t_0-2^{-2j}=k\, 2^{-2j}$ and
$$
\hat t_n:=\max\{t_n, t_\infty \}\,.
$$
Observe that the scaling property of
Brownian motion shows that there is
a constant $\rho<1$, which does not
depend on $j$ or $k$, such that
$\P[N>1]=\rho$. Moreover, the Markov property
implies that $\Pb{N\ge m+1\md N\ge m}\le \rho$
for every $m\in\N$.  In particular,
$\P[N> m]\le \rho^m$.

For every $s\in[0,\infty)$,
the map $\hf_{s}$ is measurable with respect to the
$\sigma$-field generated by
$ \xo(t)$ for $t\in[0,s]$, while $\hat t_n$ is determined by
$\xi(t)$ for $t\ge t_n$ (i.e., $\hat t_n$ is a stopping time for
the reversed time filtration).
Therefore, by the strong Markov property,
for every $n\in\N$, $s\in[t_\infty,t_0]$ and $\delta>0$,
\begin{equation*} %
\PB{|\hf_{\hat t_n}'(i\,2^{-j})|>\delta\md \hat t_n=s}
=
\PB{|\hf_{s}'(i\,2^{-j})|>\delta}\,,
\end{equation*}
and consequently,
$$
\PB{|\hf_{\hat t_n}'(i\,2^{-j})|>\delta\md \hat t_n>t_\infty}
\le
\sup_{s\in[0,1]}
\PB{|\hf_{s}'(i\,2^{-j})|>\delta}\,.
$$
This allows us to estimate,
\begin{align*}
&
\PB{\exists n\in\N\;:\; |\hf_{\hat t_n}'(i\,2^{-j})|>\delta}
\\&\qquad
\le
\Pb{|\hf_{t_\infty}'(i\,2^{-j})|>\delta} +
\sum_{n=0}^\infty \P[\hat t_n>t_\infty] \,
\PB{|\hf_{\hat t_n}'(i\,2^{-j})|>\delta\md \hat t_n>t_\infty}
\\ &\qquad
\le
 \E[N+1]\,
\sup_{s\in[0,1] }
\Pb{|\hf_s'(i\,2^{-j})|>\delta}
\le O(1)
\,\sup_{s\in[0,1]}
\Pb{|\hf_s'(i\,2^{-j})|>\delta}\,.
\end{align*}
By Corollary~\ref{derest}, and because $\sigma<(\lambda-2)/\max\{a,\lambda\}$,
this gives
\begin{equation}\label{unwant}
\PB{\exists n\in\N\;:\; \bigl|\hf_{\hat t_n}'(i\,2^{-j})\bigr|> 2^{j}\,2^{-j\sigma}/j^2}
\le
O(1)\, 2^{-2j}\,2^{-\epsilon j}\,,
\end{equation}
for some $\epsilon=\epsilon(\kappa)>0$.
Let $S$ be the rectangle
$$
S:= \bigl\{x+iy: |x|\le 2^{-j+3},\, y\in[2^{-j-1},2^{-j+3}]\bigr\}
\,.
$$
We claim that
\begin{equation}\label{subs}
H\bigl(R(j,k)\bigr)\subset
\bigcup_{n=0}^N
\hf_{\hat t_n}(S)\,,
\end{equation}
and
\begin{equation}\label{inter}
\forall n\in\N\qquad
\hf_{\hat t_n}(S)\cap\hf_{\hat t_{n+1}}(S)\ne\emptyset\,.
\end{equation}
Indeed, to prove~\eref{subs}, let $t\in [\hat t_{n+1},\hat t_n]$ and
$y\in [2^{-j-1},2^{-j}]$.
Then we may write
\begin{equation}\label{reps}
\hf_{t}(i\,y)=
\hf_{\hat t_{n+1}} %
\Bigl(g_{\hat t_{n+1}}\bigl(\hf_t(i\,y)\bigr)-\xo(\hat t_{n+1})\Bigr)
\end{equation}
We want to prove~\eref{subs} by showing that
$g_{\hat t_{n+1}}\bigl(\hf_t(i\,y)\bigr)-\xo(\hat t_{n+1})\in S$.
Set $\varphi(s):=g_{s}\bigl(\hf_t(i\,y)\bigr)$ for $s\le t$.
Then $\varphi(t)=i\,y+\xo(t)$ and by~\eref{chordal}
$$
\p_s\varphi(s)= 2\,\bigl(\varphi(s)-\xo(s)\bigr)^{-1}\,.
$$
Note that $\p_s\Im\bigl(\varphi(s)\bigr)<0$, and hence
$\Im\bigl(\varphi(s)\bigr)\ge\Im\bigl(\varphi(t)\bigr)\ge 2^{-j-1}$.
This then implies that
$|\p_s\varphi(s)|\le 2^{j+2}$, and therefore
$|t-\hat t_{n+1}|\le 2^{-2j}$ gives
$|\varphi(\hat t_{n+1})-\varphi(t)|\le 2^{2-j}$.
Since $|\xo(t)-\xo(\hat t_{n+1})|\le 2^{1-j}$,
by~\eref{reps} this gives
$ \hf_t(i\,y)\in \hf_{\hat t_{n+1}}(S) $
and verifies~\eref{subs}.  We also have~\eref{inter}, because
taking $t=\hat t_{n}$ in the above gives
$ \hf_{\hat t_{n}}(i\,y)\in \hf_{\hat t_{n+1}}(S)$.
Since $|\hf_{t}'(z)|/|\hf_{t}'(i\,2^{-j})|$ is bounded
by some constant if $z\in S$ (this follows from the Koebe distortion
theorem, see~\cite[\S1.3]{\BPomm}), we find that
$$
\diam\bigl(\hf_t(S)\bigr)\le O(1)\, 2^{-j} \,|\hf_t'(i\,2^{-j})|\,.
$$
Therefore, the relations~\eref{subs} and~\eref{inter} give
\begin{equation}\label{djk}
\begin{aligned}
d(j,k)
&
\le
\sum_{n=0}^N \diam\bigl(\hf_{\hat t_n}(S)\bigr)
\le
O(1)\, 2^{-j}\,\sum_{n=0}^N \bigl|\hf_{\hat t_n}'(i\,2^{-j})\bigr|
\\&
\le O(1)\, 2^{-j}\,N\,\max\Bigl\{|\hf_{\hat t_n}'(i\,2^{-j})\bigr|:
n=0,1,\dots,N\Bigr\}\,.
\end{aligned}
\end{equation}
By~\eref{unwant}, we get
\begin{align*}
\Pb{d(j,k)>2^{-j\sigma}}
&
\le \Ps{N>j^2}+
O(1)\, 2^{-2j}\,2^{-\epsilon j}
\\&
\le \rho^{j^2}+
O(1)\, 2^{-2j}\,2^{-\epsilon j}
\le O(1)\, 2^{-2j}\,2^{-\epsilon j}\,,
\end{align*}
which proves~\eref{Asm}.

{}From~\eref{Asm} we conclude that a.s.\ there are at most
finitely many pairs $j,k\in\N$ with $k\le 2^{2j}-1$
such that $d(j,k)> 2^{-j\sigma}$.
Hence $d(j,k)\leq C(\omega)\, 2^{-j\sigma}$ for all $j,k$,
where $C=C(\omega)$ is random (and the notation $C(\omega)$ is
meant to suggest that).
Let $(y',t')$ and $(y'',t'')$ be points in $(0,1)^2$.
Let $j_1$ be the smallest integer larger than
$\min\bigl\{-\log_2 y',-\log_2 y'',-\frac12 \log_2 |t'-t''|\bigr\}$.
Note that a rectangle
$R(j_1,k')$ that intersects the line $t=t'$ is adjacent
to a rectangle $R(j_1,k'')$ that intersects the line $t=t''$.
Consequently, $\bigl|H(y',t')-H(y'',t'')\bigr|\le
\sum_{j\ge j_1} (d(j,k'_j)+d(j,k''_j))\le O(1)\, C(\omega)\,2^{-\sigma j_1}$,
where $R(j,k'_j)$ is a rectangle meeting the line $t=t'$
and $R(j,k''_j)$ is a rectangle meeting the line $t=t''$.
This shows that for every $t_0\in[0,1)$ the limit of $H(y,t)$ as $(y,t)\to (0,t_0)$
exists, and thereby extends the definition of $H$ to $[0,\infty)\times [0,1)$.
Since $H$ is clearly continuous in $(0,\infty)\times [0,t)$,
the proof is now complete.
\QED

For $\kappa=8$, we are unable to prove Theorem~\ref{extend}.
However, a weaker result does follow easily, namely,
a.s.\ for almost every $t\ge0$ the limit $\lim_{y\searrow 0} \hf_t(i\,y)$
exists.

\begin{update}
It follows from~\cite{\BLSWlesl} and the
results of the current paper that the Theorem holds also when $\kappa=8$.
\end{update}

\section{Reduction}\label{Reduction}

The following theorem provides a criterion for hulls to be generated by
a continuous path.
In this section we do not assume that $\xo$ is
a (time scaled) Brownian motion.

\begin{theorem}\label{reduction} Let $\xo:[0,\infty)\to\R$ be continuous,
and let $g_t$ be the corresponding solution of~\eref{chordal}.
Assume that $\beta(t):=\lim_{y\searrow0} g_t^{-1}(\xo(t)+iy)$
exists for all $t\in[0,\infty)$ and is continuous.
Then $g_t^{-1}$ extends continuously to $\overline{\H}$
and $H_t$ is the unbounded connected
component of $\H\setminus \beta([0,t])$, for every $t\in[0,\infty)$.
\end{theorem}

In the proof, the following basic properties of conformal
maps will be needed.
Suppose that $g:D\to\U$ is a conformal
homeomorphism.  If $\alpha:[0,1)\to D$ is a path
such that the limit $l_1=\lim_{t\nearrow1}\alpha(t)$
exists, then $l_2=\lim_{t\nearrow1}g\circ\alpha$ exists, too.
(However, if $\alpha:[0,1)\to\U$ is a path such that
the above limit exists, it does not follow that
$\lim_{t\nearrow1}g^{-1}\circ\alpha(t)$ exists.
In other words, it is essential that the image of $g$ is a nice
domain such as $\U$.)
Moreover, $\lim_{t\nearrow1}g^{-1}(t l_2)$ exists and equals $l_1$.
Consequently, if $\tilde\alpha:[0,1)\to D$ is another path
with $\lim_{t\nearrow 1}\tilde\alpha$ existing and with
$\lim_{t\nearrow 1}g\circ\alpha(t)=\lim_{t\nearrow1}g\circ\tilde\alpha(t)$,
then $\lim_{t\nearrow 1}\alpha(t)=\lim_{t\nearrow1}\tilde\alpha(t)$.

These statements are well known and easily established, for example
with the notion of extremal length. See~\cite[Proposition 2.14]{\BPomm}
for the first statement and~\cite[Theorem 3.5]{\BAhlfors}
 implies the second claim.

\proof
Let $S(t)\subset\closure{\H}$ be the set of limit points of $g_t^{-1}(z)$
as $z\to\xo(t)$ in $\H$.
Fix $t_0\ge0$, and  let $z_0\in S({t_0})$.
We want to show that
$z_0\in\closure{\beta([0,t_0))}$,
and hence $z_0\in\beta([0,t_0])$.
Fix some $\eps>0$.  Let
$$
t':= \sup\bigl\{ t\in[0,t_0] :
{K_t}\cap\closure{ D(z_0,\eps)}=\emptyset\bigr\}\,,
$$
where $D(z_0,\eps)$ is the open disk of radius $\eps$ about $z_0$.
We first show that
\begin{equation}\label{hit}
\beta(t')\in \closure{ D(z_0,\eps)}.
\end{equation}
Indeed, $D(z_0,\eps)\cap H_{t_0}\ne\emptyset$
since $z_0\in S(t_0)$.  Let $p\in D(z_0,\eps)\cap H_{t_0}$,
and let $p'\in K_{t'}\cap \closure {D(z_0,\eps)}$.
Let $p''$ be the first point on the line segment from
$p$ to $p'$ which is in $K_{t'}$. We want to show that
$\beta(t')=p''$.
Let $L$ be the line segment
$[p,p'')$, and note that $L\subset H_{t'}$.
Hence $g_{t'}(L)$ is a curve in $\H$ terminating at a point $x\in\R$.
If $x\neq \xi(t'),$ then $g_t(L)$ terminates at points $x(t)\neq \xi(t)$
for all $t<t'$ sufficiently close to $t'$. Because $g_\tau(p'')$ has to hit
the singularity $\xo(\tau)$ at some time $\tau\leq t',$ this implies
$p''\in {K_t}$ for $t<t'$ close to $t'$.
This contradicts the definition of $t'$ and shows $x= \xi(t')$.
Now $\beta(t')=p''$ follows because the conformal map $g_{t'}^{-1}$ of $\H$
cannot have two different limits along two arcs with the same terminal point.

Having established~\eref{hit}, since $\eps>0$ was arbitrary,
we conclude that $z_0\in\closure{\beta([0,t_0))}$ and
hence $z_0\in\beta([0,t_0])$.  This gives
$S(t)\subset \beta([0,t])$ for all $t\ge0$.
We now show  that $H_t$ is the unbounded component of
$\H\setminus \closure{ \bigcup_{\tau\leq t} S(\tau)}$.
First, $H_t$ is connected and
disjoint from $\closure{\bigcup_{\tau\leq t} S(\tau)}$.
On the other hand, as the argument in the previous
paragraph shows, $\p H_t\cap\H$ is contained in
$\closure{\bigcup_{\tau\leq t} S(\tau)}$.  Therefore,
$H_t$ is a connected component of
$\H\setminus \closure{\bigcup_{\tau\leq t} S(\tau)}$; that is,
$H_t$ is the unbounded connected component of
$\H\setminus \beta([0,t])$.  Since $\beta$ is a continuous
path, it follows from \cite[Theorem 2.1]{\BPomm} that $g_t^{-1}$ extends
continuously to $\closure{\H}$ (which also proves that $S(t)=\{\beta(t)\}$).
\QED

\section{Continuity}\label{scont}

We have now established all the results needed to show that the \SLEk/ trace is
a continuous path a.s.

\begin{theorem}[Continuity]\label{slecont} 
Let $\kappa\neq 8$.
The following statements hold almost surely.
For every $t\ge 0$ the limit
$$
\gamma(t):=\lim_{z\to 0,\, z\in\H}\hf_t(z)
$$
exists, 
 $\gamma:[0,\infty)\to\closure{\H}$ is a continuous path,
and $H_t$ is the unbounded component of $\H\setminus\gamma\bigl([0,t]\bigr)$.
\end{theorem}

We believe the theorem to be valid also for $\kappa=8$.
(This is stated as Conjecture~\ref{slecont8}.)
Despite repeated efforts, the proof eluded us.

\begin{update}
This extension to $\kappa=8$ is proved in~\cite{\BLSWlesl}.
\end{update}

\proofof{Theorem \ref{slecont}}
By Theorem~\ref{extend}, a.s.\ $\lim_{y\searrow 0}\hf_t(iy)$
exists for all $t$ and is continuous.
Therefore we can apply Theorem~\ref{reduction}, and the theorem follows.
\QED

It follows from Theorem~\ref{slecont} that $f_t$ extends continuously to $\overline\H$ a.s.
The next result
gives more information about the regularity of $f_t$ on $\overline\H$. It neither
follows from Theorem~\ref{slecont}, nor does it imply \ref{slecont}.

\begin{theorem}[H\"older Continuity]\label{holdercont}
For every $\kappa\ne4$ there is some $h(\kappa)>0$ such that
for every bounded set $A\subset\H$ and every $t>0,$ a.s.\
$\hf_t$ is H\"older continuous with exponent $h(\kappa)$ on $A,$
$$|\hf_t(z)-\hf_t(z')|\leq C |z-z'|^{h(\kappa)}\qquad\forall\quad z,z'\in A,$$
where $C=C(\omega,t,A)$ is random and may depend on $t$ and $A$.
Moreover, $\lim_{\kappa\searrow0} h(\kappa)=\frac12$ and
$\lim_{\kappa\nearrow\infty} h(\kappa)=1$.
\end{theorem}

Since $f_t(z)-z\to0$ as $z\to\infty,$ it easily follows that for every $t$ a.s.
\begin{equation}
|\hf_t(z)-\hf_t(z')|\leq C(\omega,t)\max(|z-z'|, |z-z'|^{h(\kappa)}).\label{globalholder}
\end{equation}

\medskip
We don't believe that the theorem holds for $\kappa=4$, for then
the trace is a simple path which ``almost'' touches itself.

\begin{update}
For $\kappa\leq4$, the fact that $\gamma$ is a simple path in $\H$, 
Theorem \ref{simplepath} below, implies $h(\kappa)\leq 1/2$. Thus the
estimate $\lim_{\kappa\searrow0} h(\kappa)=\frac12$ is best possible.
On the other hand, this non-smoothness of $\hf_t$ is localized at
$\hf_t^{-1}(0)$: Joan Lind (manuscript in preparation) has shown that
the H\"older exponent $h(\kappa)$ of $(\hf_t(\sqrt z))^2$ satisfies
$\lim_{\kappa\searrow0} h(\kappa)=1$.
\end{update}

\proof
Fix $\kappa\ne 4$ and $t>0$. By scaling, we may assume $0<t\leq1$
and $A=[-1,1]\times (0,1]$. 
Denote
$$z_{j,n}=(j+i) 2^{-n},\quad 0\leq n <\infty,\quad -2^{n}< j <2^{n}.$$
We will first show that there is an $h=h(\kappa)>0$ such that a.s.
\begin{equation}\label{derivativegrowth}
|\hf_t'(z_{j,n})|\leq C(\omega,t) 2^{n (1-h)},\qquad \forall j,n.
\end{equation}
Using Corollary~\ref{derest} with $\delta=2^{-n h}$ we have
$$\Pb{ |\hf_t'(z_{j,n})|\ge 2^{n (1-h)} } \le
C(\kappa,b)\, (1+2^{2\,n})^b\, 2^{-n(1-h)\lambda}\,\hh(2^{-n h},a-\lambda).$$
Hence
$$\sum_{n=0}^{\infty} \sum_{j=-2^{n}}^{2^{n}}
\Pb{ |\hf_t'(z_{j,n})|\ge 2^{n (1-h)} } < \infty$$
provided that
$$1+2b-(1-h)\lambda <0 \quad \text{and} \quad a-\lambda\le0,$$
or that
$$1+2b-\lambda+a h<0 \quad \text{and} \quad a-\lambda\ge0.$$
If $0<\kappa\le 12,$ $b=1/4+1/\kappa$ and $h<(\kappa-4)^2/((\kappa+4)(\kappa+12)),$ the first condition is satisfied.
For $\kappa>12$, $b=4/\kappa$ and $h<1/2-4/\kappa$ the second condition is
satisfied, and~\eref{derivativegrowth} follows.

To see that one can actually achieve $h(\kappa)\to 1/2$ as $\kappa\searrow 0,$
set $b:=-1/2 + \sqrt{1/2+2/\kappa}$ for $0<\kappa < 4$,
and let $h$ be smaller than but close to $(\lambda-1-2b)/\lambda$.
To get $\lim_{\kappa\nearrow\infty} h(\kappa)=1$, take
$b:=(\sqrt{2\kappa}\sqrt{\kappa^2+10\kappa+16} -2\kappa)/(\kappa^2+4\kappa)$
for $\kappa\ge 2(3+\sqrt{17})$ and let $h$ be smaller than but close
to $(\lambda-1-2b)/a$. 

{}From the Koebe distortion Theorem and~\eref{derivativegrowth} we obtain
$$|\hf_t'(z)|\le O(1)\,C(\omega,t)\, y^{h-1}$$ 
for all $z\in A$. It is well-known and easy to see, by integrating $|f'|$
over the hyperbolic geodesic from $z$ to $z'$ (similarly to the end of the
proof of Theorem~\ref{extend}), that this implies H\"older
continuity with exponent $h$ on $A$. The theorem follows.
\QED

The following corollary is an immediate consequence of Theorem~\ref{holdercont}.
In Section~\ref{dimensions} below, we will present more precise estimates.

\begin{corollary}\label{smallbd} For $\kappa\ne4$ and every $t,$
the Hausdorff dimension of $\partial K_t$ a.s.\ satisfies
$$\mathrm{dim}\, \partial K_t < 2.$$
In particular, a.s.\
$\area\, \partial K_t=0$.
\end{corollary}

\proof By \cite{\BJonMak} (see also \cite{\BKR} for an 
easier proof), the boundary of
the image of a disk under a
H\"older continuous conformal map has Hausdorff dimension bounded away from 2.
Consider the conformal map $T(z)=(z-i)/(z+i)$ from $\H$ onto
$\U$.
By~\eref{globalholder}, $T\circ f_t\circ T^{-1}$ is a.s.\
H\"older continuous in $\U$. Since $T$ preserves Hausdorff dimension,
the corollary follows.
\QED

\section{Phases}\label{sphases}

In this section, we will investigate the topological behavior of \SLE/, 
and will identify three very different phases for the parameter $\kappa$,
namely, $[0,4]$, $(4,8)$, and $[8,\infty)$.

The following result was conjectured in~\cite{\BSchSLE}. There, it was proved that
for $\kappa>4,$ a.s.\ $K_t$ is not a simple path. The proof was based on the calculation
of the harmonic measure $F(x)$ below, which we will repeat here for the convenience
of the reader.

\begin{theorem}\label{simplepath}
In the range $\kappa\in[0,4]$, the \SLEk/ trace $\gamma$
is a.s.\ a simple path and $\gamma[0,\infty)\subset\H\cup\{0\}$.
\end{theorem}

\begin{lemma}\label{Rpart}
Let $\kappa\in[0,4]$, and let $\gamma$ be the
\SLEk/ trace.
Then a.s.\
$$
\gamma[0,\infty)\subset \H\cup\{0\}\,.
$$
\end{lemma}

Set $Y_x(t):=g_t(x)-\xo(t)$, $x\in\R$, $t\ge 0$. 
Then $Y_x(t)/\sqrt\kappa$ is a Bessel process
of dimension $1+4/\kappa$.
The basic theory of Bessel processes
then tells us that a.s.\ $Y_x(t)\ne 0$ for $t\ge 0$, 
$x\ne 0$ and $\kappa\in[0,4]$, which amounts to $x\ne K_t$.
However, we will prove this for the convenience of the reader unfamiliar
with Bessel processes.

\proof
Let $b>a>0$.
For $x\in[a,b]$,
let $T=T_x:=\inf\{t\ge 0:Y_x(t)\notin(a,b)\}$.
Set
$$
\hat F(x) :=
\begin{cases}
x^{(\kappa- 4)/\kappa}& \kappa\ne 4,\\
\log(x)&\kappa=4\,,
\end{cases}
$$
and
$$
F(x):= \frac{\hat F(x)-\hat F(a)}{\hat F(b)-\hat F(a)}\,.
$$
\Ito/'s formula shows that $F(Y_x(t\wedge T))$
is a local martingale, and since $F$ is bounded in $[a,b]$,
this is a martingale.  Consequently, the optional sampling theorem
gives
$F(x)=\Eb{F(Y_x(T)}=\Pb{Y_x(T)=b}$.

Note that $F(x)\to 1$ when $a\searrow 0$.
(That's where the assumption $\kappa\le 4$ is crucial.)
Hence, given $x>0$, a.s.\ for all $b>x$, there is some $s>0$ such that
$g_t(x)$ is well defined for $t\in[0,s]$,
$\inf\bigl\{Y_x(t):t\in[0,s]\bigr\}>0$ and $Y_x(s)=b$.
Note that the \Ito/ derivative of $Y_x(t)$ (with respect to $t$) is
$$
dY_x = (2/Y_x)\,dt + d\xo\,.
$$
It follows easily
that a.s.\ $Y_x(t)$ does not escape to $\infty$ in finite time.
Observe that if $x'>x$, then $Y_{x'}(t)\ge Y_x(t)$.
Therefore, a.s.\ for every $x>0$
we have $Y_x(t)$ well defined and in $(0,\infty)$ for all $t\ge 0$.
This implies that a.s.\ for every $x>0$ and every $s>0$
there is some neighborhood $N$ of
$x$ in $\C$ such that the differential equation~\eref{chordal}
has a solution in the range $z\in N$, $t\in[0,s]$.
This proves that a.s.\ $\gamma[0,\infty)$
does not intersect $(0,\infty)$.
The proof that it a.s.\ does not intersect $(-\infty,0)$ is the
same.
\QED

\proofof{Theorem \ref{simplepath}}
Let $t_2>t_1>0$.  The theorem will be established by proving that
$\gamma[0,t_1]\cap\gamma[t_2,\infty)=\emptyset$.
Let $s\in (t_1,t_2)$ be rational, and set
$$
\hat g_t(z):= g_{t+s}\circ g_{s}^{-1}\bigl(z+\xo({s})\bigr) - \xo(s).
$$
By Proposition~\ref{funda} $(\hat g_t:t\ge 0)$ has the same distribution
as $(g_t:t\ge 0)$.
Let $\hat\gamma_s (t)$ be the trace for the collection
$(\hat g_t:t\ge 0)$; that is,
$$
\hat\gamma_s(t)= g_{s} \circ g_{t+s}^{-1}\bigl(\xo(t+s)\bigr)-\xo(s)
=g_s\circ \gamma(t+s)-\xo(s)\,.
$$
We know from Lemma~\ref{Rpart}
that a.s.\ for all rational $s\in[t_1,t_2]$,
$$
\hat\gamma_s[0,\infty)\subset \H\cup\{0\}\,.
$$
Consequently, by applying $g_s^{-1}$, we conclude
that a.s.\ for every rational $s>0$,
$$
\gamma[s,\infty)\cap (\R\cup K_s) = \{\gamma(s)\}\,.
$$
Since $g_{t_2}\circ g_{t_1}^{-1}$ is
not the identity, it follows from
Theorem~\ref{slecont} that there is some
$s\in(t_1,t_2)$ such that
$\gamma(s)\notin\R\cup K_{t_1}$,
and this holds for an open set of $s\in (t_1,t_2)$.
In particular, there is some rational $s\in(t_1,t_2)$
with $\gamma(s)\notin\R\cup K_{t_1}$.
Since $\gamma[0,t_1]\subset K_{t_1}\cup\R$,
it follows that
$\gamma[0,t_1]\cap\gamma[t_2,\infty)=\emptyset$,
proving the theorem.
\QED


Recall the definition of the hypergeometric function in $|z|<1$,
\begin{equation}\label{hypergeom}
{}_2F_1(\eta_0,\eta_1,\eta_2,z)
:=\sum_{n=0}^\infty
\frac{(\eta_0)_n\,(\eta_1)_n}{(\eta_2)_n\,n!}\,z^n\,,
\end{equation}
where $(\eta)_n:= \prod_{j=0}^{n-1} (\eta+j)$.
Let
$$
\hG(x+i\,y)=\hG_{a,\kappa}(x+i\,y):=
{}_2F_1\bigl(\eta_0,\eta_1,1/2,x^2/(x^2+y^2)\bigr)\,,
$$
where
$$
\eta_j=\eta_j(a,\kappa):=
\frac 12-\frac2\kappa-
(-1)^j\,\frac{\sqrt{32\,a\,\kappa+(2\,\kappa-8)^2}}{4\,\kappa}
\,.
$$

\label{Showclose}
\begin{lemma}\label{howclose}
Let $\kappa> 0$ and $z=x+iy\in\H$.  Then the limit
$$
Z(z):=\lim_{t\nearrow \tau(z)}
 y\,\bigl|g_t'(z)\bigr|/\Im\bigl( g_t(z)\bigr)
$$
exists a.s.
We have $Z(z)=\infty$ a.s.\ if $\kappa\ge 8$
and $Z(z)<\infty$ a.s.\ if $\kappa<8$.
Moreover, for all $a\in\R$  we have
\begin{equation*}
\begin{aligned}
\Eb{Z(z)^a}
=
\begin{cases}
    \hG_{a,\kappa}(z)/\hG_{a,\kappa}(1)& a<1-\kappa/8,\ \kappa<8,\\
    \infty& a\ge 1-\kappa/8,\ a>0,\\
    0& a<0,\ \kappa\ge8.
\end{cases}
\end{aligned}
\end{equation*}
\end{lemma}

Before proving the lemma, let's discuss its geometric meaning.
For $z_0\in\H$ and $t<\tau(z_0)$ we have $z_0\in H_t$.
Set $r_t:=\dist(z_0,\partial H_t)$, and let
$\phi:\H\to\U$ be a conformal homeomorphism satisfying $\phi\bigl(g_t(z_0)\bigr)=0$.
Note that $\Bigl|\phi'\big(g_t(z_0)\bigr)\Bigr|=\bigl(2 \,\Im\, g_t(z_0)\bigr)^{-1}$.
Hence, the Schwarz Lemma applied to the map
$z\mapsto\phi\circ g_t(z_0+r_t\, z)$ proves
$r_t\,\bigl|g_t'(z_0)\bigr|\le 2\,\Im\, g_t(z_0)$.
On the other hand, the Koebe 1/4 Theorem applied to the function
$g_t^{-1}\circ \phi^{-1}$ implies that
$\Im\, g_t(z_0) \le 2 \, r_t\, |g'_t(z_0)|$.
%
Therefore, since $\lim_{t\nearrow\tau(z_0)} r_t 
=\dist\bigl(z_0,\gamma[0,\infty)\cup\R\bigr)$,
\begin{equation}\label{dtz}
 Z(z_0)^{-1}\,\Im\, z_0/2\le
\dist\bigl(z_0,\gamma[0,\infty)\cup\R\bigr)\le 2 \,Z(z_0)^{-1}\,\Im\, z_0\,.
\end{equation}
In particular, the lemma tells us that $\gamma[0,\infty)$ is dense in
$\H$ iff $\kappa\ge 8$.
(For $\kappa=8$ we have not proven that $\gamma$ is a path.
In that case, $\gamma[0,\infty)$ is defined as the union
over $t$ of all limit points of sequences $\hf_t(z_j)$ where
$z_j\to0$ in $\H$.  By an argument in the proof of
Theorem~\ref{reduction}, 
$\closure{\bigcup_{t\ge 0}\p K_t}=\closure{\gamma[0,\infty)}$.)

The way one finds the function $\hG$ is explained within the proof below, and
is similar to the situation in Theorem~\ref{derexp}.

\proof
Let $G(z)=G_{a,\kappa}(z):=\Eb{Z(z)^a}$.
Fix some $\iz=\ix+i\,\iy\in\H$, and
abbreviate $Z:=Z(\iz)$.
Let $z_t=x_t+i\,y_t:=g_t(\iz)-\xo(t)$.
{}From~\eref{chng} we find that
$$
\p_t \log \bigl|g_t'(z)/y_t\bigr|= 4\, y_t^2\,|z_t|^{-4}\,,
$$
which implies that
\begin{equation}\label{Zis}
Z = \exp\Bigl(\int_{0}^{\tau(\iz)}
4\,y_t^2\, |z_t|^{-4}\,dt\Bigr)
 =  \exp\Bigl(\int_{0}^{\tau(\iz)}
\frac {-2\,y_t^2}{x_t^2+y_t^2}\,d\log y_t\Bigr)\,.
\end{equation}
In particular, the limit in the definition of $Z$ exists
a.s.\ and $G(\iz)=\E[Z^a]$ (possibly $+\infty$).

Let $\ev F_t$ denote the $\sigma$-field generated by
$\bigl(\xo(s):s\le t\bigr)$.
A direct application of \Ito/'s formula shows that
$$
M_t:= \bigl(\iy\,|g_t'(\iz)|/y_t\bigr)^a\,\hG(z_t)
$$
is a local martingale\footnote{
The Markov property can be used to show that
$\bigl(\iy\,|g_t'(\iz)|/y_t\bigr)^a\,G(z_t)$ is a local martingale.
One then assumes that $G$ is smooth and applies \Ito/'s formula, to obtain
 a PDE satisfied by $G$.  Scale invariance can then be used to transform
the PDE to the ODE~\eref{G} appearing below. 
Then one obtains $\hG$ as a solution.  That's how
the function $\hG$ was arrived at.}.

As before, let $u_t=\log y_t$.
Note that $w(u)=w_t:=x_t/y_t$ is a time-homogeneous diffusion process
as a function of $u$.  It is easy to see that $\tau(\iz)$ is a.s.\ the time
at which $u$ hits $-\infty$. (For example, see the argument following~\eref{Tdef}.)

If $t_1<\tau(\iz)$ is a stopping time, then by~\eref{Zis} we may write
\begin{equation}\label{splitZ}
\log Z(\iz)=
 \int_{0}^{t_1} + \int_{t_1}^{\tau(\iz)}
4\,y_t^2\, |z_t|^{-4}\,dt\,.
\end{equation}
Note that if $s>0$ then the distribution of $Z(\iz)$ is the
same as that of $Z(s\,\iz)$, by scale invariance.
It is clear that for every $z'\in\H$ there is positive probability
that there will be some first time ${t_1}<\tau(\iz)$ such that
$z_{t_1}/|z_{t_1}|=z'/|z'|$.
Conditioned on $t_1<\tau(\iz)$ and on $\ev F_{t_1}$,
the second integral in~\eref{splitZ}
has the same distribution as $\log Z(z')$, by the strong Markov property
and scale invariance.  Consequently, we find that $\Pb{Z(\iz)>s}\ge c \,\Pb{Z(z')>s}$
holds for every $s$, where $c=c(\iz,z',\kappa)>0$.
In fact, later we will need the following stronger statement
(which is not needed for the present lemma, but we find it convenient to
give the argument here). 
Let $\overline w\in(0,\infty)$.  Then
assuming $|w_0|\le \overline w$, there is a $c=c(\overline w,\kappa)>0$
such that
\begin{equation}\label{comp}
\forall s>0\,,\qquad \Pb{Z(i)>s}\ge \Pb{Z(\iz)>s}\ge c \,\Pb{Z(i)>2\,s}\,.
\end{equation}
The right hand inequality is proved by taking $t_1$ to be the first time
smaller than $\tau(z)$ such that
$\Re\, z_{t_1}=0$ and  the first integral in~\eref{splitZ} is
greater than $\log 2$, on the positive probability event that there is such a time.
The left hand inequality is clear, for if $|\Re\, g_t(i)/\Im\, g_t(i)|$
never hits $|w_0|$, then we have $Z(i)=\infty$, by~\eref{Zis}.

Note that $\hG(i)=1$.  Assume, for now, that $G(i)<\infty$.
We claim that in this case $\hG(x+i)>0$ for each $x\in\R$.
Let $s_0>0$, let $T:=\inf\{t\ge 0: |w_t|=s_0\}$.
Since $w(u)$ is time-homogeneous, clearly, $T<\tau(\iz)$.
If $\hG(x+i)\le 0$ for some $x\in\R$,
we may choose $s_0:=\inf\{x>0:\hG(x+i)=0\}$.
Then $M_{T}=0$ if $|w_0|\le s_0$.
Note that 
$1\le\iy\,|g_{T}'(\iz)|/y_{T}<Z$.
Since $M_t$ is a local martingale on $t< \tau(\iz)$, there is an 
increasing sequence
of stopping times $t_n<\tau(\iz)$ with $\lim_n t_n=\tau(\iz)$ a.s.\ 
such that $M_{t\wedge t_n}$ is a martingale.
Set $\bar t_n:=T\wedge t_n$.  The optional
sampling theorem then gives for $\iz=i$
\begin{equation*}
\begin{split}
1=\hG(\iz)=M_0=\Es{M_{\bar t_n}}=
\Eb{\bigl(\iy\,|g_{\bar t_n}'(\iz)|/y_{\bar t_n}\bigr)^a\,\hG(z_{\bar t_n})}
\qquad
\\
\le
\Eb{\max\{1, Z^a\}\,\hG(z_{\bar t_n})}\,.
\end{split}
\end{equation*}
Since $\hG(z_{t\wedge T})$ is bounded, $\lim_{n\to\infty}\hG(z_{\bar t_n})= 0$ a.s.,
and $\E[Z^a]<\infty$, we have
$\lim_{n\to\infty}\Eb{\max\{1, Z^a\}\,\hG(z_{\bar t_n})}=0$,
contradicting the above inequality.
  Consequently, $\hG(x+i)>0$ for every $x\in\R$.

We maintain the assumption $G(i)<\infty$.
By~\eref{comp}, this implies $G(z)<\infty$ for all $z\in\H$. 
Let $s_1,s_2,\dots$ be a sequence in $(0,\infty)$ with $s_n\to\infty$ such that
$L:=\lim_n\hG(s_n+i)$ exists (allowing $L=\infty$).
Now set $T_m:=\inf\{t\ge 0: |w_t|=s_m\}$. 
The argument of the previous paragraph shows that
$M_0=\Eb{M_{T_m\wedge t_n}}\to \Eb{M_{T_m}}$,
as $n\to\infty$.  Therefore,
$$
\hG(\iz)=M_0=\Eb{M_{T_m}}\to L\,G(\iz)\,,\qquad m\to\infty\,,
$$
unless $L=\infty$, which happens if and only if $G(\iz)=0$.
Moreover, $L$ does not depend on the particular sequence $s_n$, and
we may write $\hG(1)$ instead of $L$.
Thus, when $G(i)<\infty$,
\begin{equation}\label{firstline}
G(\iz)=\hG(\iz)/\hG(1)\,.
\end{equation}

When $a=1-\kappa/8$, $\hG$ simplifies to $(y/|z|)^{(8-\kappa)/\kappa}$,
with $\hG(1)=\infty$ if $\kappa>8$ and $\hG(1)=0$ if $\kappa<8$.
It follows that $G(z)=\infty$ if $a=1-\kappa/8$ and $\kappa<8$.
Clearly, $G(z)=\infty$ also if $a>1-\kappa/8$, and $\kappa<8$.
Similarly, we have $G(z)=0$ if $\kappa>8$ and $a=1-\kappa/8$.
This implies that $Z=\infty$ a.s.\ when $\kappa>8$, which implies
$G=0$ when $\kappa>8$ and $a<0$.
There are various ways to argue that $Z=\infty$ a.s.\ when $\kappa=8$.
In particular, it is easy to see that $\hG_{\frac14,8}(1)=\infty$ 
(if $a_n$ denotes the $n$-th coefficient of 
${}_2F_1(\frac14,\frac14,\frac12,z)$, then 
$(n+1) a_{n+1} = n a_n b_n$ with $b_n=(n+1/4)^2/(n(n+1/2))$ and it follows 
from $\prod_n b_n >0$ that $n a_n$ is bounded away from zero) 
so that the above argument applies.
It remains to prove
that $G(i)<\infty$ when $a<1-\kappa/8$ and $\kappa<8$.

Assume now that $\kappa<8$ and $a<1-\kappa/8$
and $a$ is sufficiently close to $1-\kappa/8$ so that
$\eta_1>0$.
We prove below that in this range
\begin{equation}\label{lbd}
\inf_{z\in\H}\hG(z)>0\,.
\end{equation}
Since $\hG(\iz)=M_0=\Eb{M_{t_n}}$
when $\iz=i$, we immediately get $G(i)<\infty$ from~\eref{lbd}.
Since $G(i)<\infty$ for some $a$ implies $G(i)<\infty$ for any smaller $a$,
it only remains to establish~\eref{lbd} in the above specified range.

To prove~\eref{lbd},
first note that $\eta_0\le\eta_1$, $0<\eta_1<1/2$ and
$\eta_0+\eta_1= 1-4/\kappa<1/2$.
By~\cite[2.1.3.(14)]{\BhigherFunc}, when $\eta_2>\eta_1+\eta_0$
and $\eta_2>\eta_1>0$, the right hand side
of~\eref{hypergeom} converges at $z=1$ and is equal to
\begin{equation}\label{gauss}
\frac{\Gamma(\eta_2)\,\Gamma(\eta_2-\eta_1-\eta_0)}
{\Gamma(\eta_2-\eta_1)\,
\Gamma(\eta_2-\eta_0)}\,.
\end{equation}
This gives $\hG(1)>0$.
Now write
$$
g(x):= \hG_{a,\kappa}(x+i)/\hG_{1-\kappa/8,\kappa}(x+i)
=(1+x^2)^{(8-\kappa)/(2\kappa)}\,\hG_{a,\kappa}(x+i)\,.
$$
Since $\hG$ satisfies
\begin{equation}\label{G}
\frac{4\, a\, y^2}{|z|^4}\,\hG+\frac{\kappa}{2}\,\p_x^2 \hG
+ \frac{4 x}{|z|^2}\,\p_x \hG=0
\end{equation}
for $y=1$,
it follows that
\begin{equation}\label{g}
(\kappa-8+8a)\,g+2\,(\kappa-4)\,x\,(1+x^2)\,g'+\kappa\,(1+x^2)^2\,g''=0\,.
\end{equation}
It is immediate that $g(0)=1$, $g'(0)=0$ and $g''(0)>0$,
e.g., by differentiation or by considering~\eref{G}.
Consequently, $g>1$ for small $|x|>0$.  Suppose that
there is some $x\ne0$ such that $g(x)=1$.  Then there must
be some $x_0\ne0$ such that $g(x_0)>1$ and $g$ has a local maximum at
$x_0$.  Then $g'(x_0)=0\ge g''(x_0)$.  But this contradicts~\eref{g},
giving $g>1$ for $x\ne 0$. Together with $\hG(1)>0$,
this implies~\eref{lbd}, and completes the proof.
\QED

We say that $z\in\closure\H$ is {\bf swallowed} by the
\SLE/ if $z\notin \gamma[0,\infty)$ but
$z\in K_t$ for some $t\in(0,\infty)$.  The smallest such
$t$ will be called the {\bf swallowing time} of $z$.

\begin{theorem}\label{tswal}
If $\kappa\in(4,8)$ and $z\in\closure{\H}\setminus\{0\}$,
then a.s.\ $z$ is swallowed.
If $\kappa\notin(4,8)$, then a.s.\ no $z\in\closure\H\setminus\{0\}$
is swallowed.
\end{theorem}

\begin{lemma}\label{posvol}
Let $z\in\closure{\H}\setminus\{0\}$.
If $\kappa>4$, then $ \Pb{\tau(z)<\infty}=1 $.
\end{lemma}

In fact, in this range $\kappa>4$,
the stronger statement that a.s.\ for all $z\in\closure{\H}$ we have
$\tau(z)<\infty$ holds as well.  This follows
from Theorem~\ref{transient} (and its extension to $\kappa=8$ based
on~\cite{\BLSWlesl}).

\proof
Let $b:=1-4/\kappa$, and $\theta:= \exp\bigl(i \pi (1-b)/2\bigr)$.
Set
$$
h(z):=\Im\bigl(\theta\,z^b\bigr), \qquad z\in \closure{\H}\,.
$$
A direct calculation shows that the drift term in \Ito/'s formula for
$\bigr(g_t(z)-\xo(t)\bigr)^b$ is zero.
Consequently, $h\bigl(g_t(z)-\xo(t)\bigr)$ is a local martingale.
Note that $\kappa>4$ implies that $h(z)\to\infty$ as
$z\to\infty$ in $\closure{\H}$, that $h(z)>0$ in $\closure{\H}\setminus\{0\}$
and $h(0)=0$.

Fix some $z_0\in\H$ and some $\epsilon>0$.
For $t\ge 0$ set $z_t:=g_t(z_0)-\xo(t)$.
We shall prove that there is some $T=T(z_0,\epsilon)$ such
that $\Pb{z_0\in K_T}>1-\epsilon$.
Indeed, let $R>1$ be sufficiently large so that
$\epsilon\,\inf\{h(z):|z|=R\}> 2h(z_0)$.
Let
$$
\tau_R:=\inf\Bigl(\bigl\{t\in[0,\tau(z_0)):|z_t|\ge R\bigr\}\cup
\bigl\{\tau(z_0)\bigr\}\Bigr)\,.
$$

We claim that if $T=T(z_0,\epsilon)$ is sufficiently
large, then $\Pb{\tau_R>T}<\epsilon/2$.
Let $t\ge 0$ and let $\ceil t:=\inf \{k\in \Z:k\ge t\}$.
Say that $t$ is a {\bf sprint time} if
$|\xo(t)-\xo(\ceil t)|=R$ and there is no
$s\in[t,\ceil t]$ such that $|\xo(s)-\xo(\ceil t)|=R$.
Note that for every integer $k\ge 0$ there is positive
probability that there is a sprint time in $[k,k+1]$
and this event is independent from $\bigl(\xo(t):t\in[0,k]\bigr)$.
Therefore, there are infinitely many sprint times a.s.
Let $ t_0<t_1<t_2<\cdots$ be the sequence of sprint times.
Let $\ev A_j$ be the event  that $\xo(\ceil {t_j})>\xo(t_j)$.
Observe that for all $j\in\N$,
\begin{equation}\label{unbia}
\Pb{\ev A_j\md t_j,\ \ev F_{t_j}}=1/2\,,
\end{equation}
where $\ev F_{t}$ denotes the $\sigma$-field generated by
$\bigl(\xo(s):s\in[0,t]\bigr)$, because
given $t_j$ and $\bigl(\xo(t):t\in[0,t_j]\bigr)$,
the distribution of $\bigl(\xo(t):t>t_j\bigr)$ is
the same as that of $\bigl(2\xo(t_j)-\xo(t):t>t_j\bigr)$,
by reflection symmetry.
Let $x_t:=\Re (z_t)$, and note that
$$
dx_t = 2 x_t\,|z_t|^{-2}\,dt - d\xo(t)\,.
$$
Consequently,  for $s\ge t_j$,
\begin{equation}\label{gogo}
x_s=x_{t_j}+
2\int_{t_j}^s x_t\,|z_t|^{-2}\,dt-\xo(s)+\xo(t_j)\,.
\end{equation}
On  $\ev A_j$ we have $\xo(s)-\xo(t_j)>0$ for all $s\in[t_j,\ceil{t_j}]$.
Hence, it follows from~\eref{gogo}
that on $\ev A_j$ and assuming $x_{t_j}\le 0$,
we have $x_s\le 0$ for all $s\in[t_j,\ceil{t_j}]$ and moreover
$x_{\ceil{t_j}}<-R$.
Similarly, if $x_{t_j}\ge 0$ and
$\neg\ev A_j$ holds, then $x_{\ceil{t_j}}>R$.
Hence, by~\eref{unbia} we have
$\Pb{\tau_R>\ceil{t_j}}\le 2^{-j}$, which implies
that we may choose $T$ such that $\Ps{\tau_R>T}<\epsilon/2$.
Set $\tau^*:=\tau_R\wedge T$, where $T$ is as above.

Since $h(z_t)$ is a local martingale, and $h(z_{t\wedge\tau^*})$ is
bounded, $h(z_{t\wedge\tau^*})$ is a martingale.
Therefore,
$$
h(z_0)=\E[h(z_{\tau^*})]\ge \inf\{h(z):|z|=R\}\,\Pb{|z_{{\tau^*}}|=R}\,.
$$
By our choice of $R$, this gives
$\Pb{|z_{\tau^*}|=R}<\epsilon/2$.  Since $\Ps{\tau_R>T}<\epsilon/2$,
it follows that $\Ps{\tau(z_0)\le T}\ge 1-\epsilon$, and the lemma
is thereby established.
\QED

We now prove a variant of Cardy's formula.  (See~\cite{\BCardySurvey}
for a survey of Cardy's formula.)
A proof of (a generalized) Cardy's formula for \SLEk/ appears in~\cite{\BLSWi}.
However the following variant appears to be different, except in
the special case $\kappa=6$, in which they are equivalent
due to the locality property of \SLEkk6/.
(See~\cite{\BLSWi} for an explanation of the locality property.)

\begin{lemma}\label{ncard}
Fix $\kappa\in(4,8)$ and
let
$$
X:=\inf\bigl([1,\infty)\cap\gamma[0,\infty)\bigr)
\,.
$$
Then  for all $s\ge 1$,
\begin{equation}\label{card}
\P[X\ge s]=
\frac{
4^{(\kappa-4)/\kappa}\,\sqrt{\pi}\,
{}_2F_1(1-4/\kappa,2-8/\kappa,2-4/\kappa,1/s) \,s^{(4-\kappa)/\kappa}
}{\Gamma(2-4/\kappa)\,\Gamma(4/\kappa-1/2)}
\,.
\end{equation}
\end{lemma}

The proof will be brief, since it is similar to
the detailed proof of the more elaborate result in~\cite{\BLSWi}.

\proof
Let $s>1$, and set
$$
Y_t:=\frac{ g_t(1)-\xo(t)}{g_t(s)-\xo(t)}
\,,
$$
for $t<\tau(1)$.
It is not hard to see that
$\lim_{t\nearrow\tau(1)} Y_t=1$ if $X>s$ and
$\lim_{t\nearrow\tau(1)} Y_t=0$ if $X<s$:
Just notice that for $X>s$ and $t\nearrow\tau(1),$ the extremal
distance within $H_t$ from $[1,s]$ to $[\Re\,\gamma(t),\gamma(t)]$ tends
to $\infty$, whereas for $X<s$ we have $\tau(1)<\tau(s)$.
Direct inspection shows that the right hand side of~\eref{card}
is zero when $s=\infty$ and is $1$ when $s=1$.
In particular, it is bounded for $s\in[1,\infty)$.
Moreover, \Ito/'s formula shows that when we plug $1/Y_t$
in place of $s$ in the right hand side of~\eref{card},
we get a local martingale.  Since $Y_0=1/s$ and $Y_t\in(0,1)$
for $t<\tau(1)$, it therefore suffices
to prove that $\lim_{t\nearrow\tau(1)}Y_t\in\{0,1\}$ a.s.
The only way this could fail is if $X=s$ occurs with
positive probability.   However, scale invariance
shows that $\P[X=s']\ge\P[X=s]$ when $s'\in(1,s)$.
Hence $\P[X=s]=0$ for all $s\ne 1$.
\QED

\proofof{Theorem~\ref{tswal}}
The case $\kappa\in [0,4]$ is covered by Theorem~\ref{simplepath}
and Theorem~\ref{slecont}.

Suppose that $\kappa\in(4,8)$ and $z\in\closure{\H}$.
Lemma~\ref{posvol} shows that
$z\in\bigcup_{t>0}K_t$ a.s.
If $z\in\H$, then Lemma~\ref{howclose} and~\eref{dtz} show
a.s.\ $z\notin\gamma[0,\infty)$.
 Lemma~\ref{ncard} shows that $1\notin\gamma[0,\infty)$
a.s., and the same follows for every $z\in\R\setminus\{0\}$ by
scale and reflection invariance.  This covers the range $\kappa\in(4,8)$.

Now suppose that $z_0$ is swallowed, and let $t_0$ be the swallowing time.
If $\delta>0$ is smaller than the distance from $z_0$ to $\gamma[0,t_0]$,
then all the points within distance $\delta$ from $z_0$ are swallowed
with $z_0$.
Therefore, if the probability that some $z_0\in\closure\H\setminus\{0\}$ is
swallowed is positive, then by~\eref{dtz} there is some $z_0$ such that
$\Pb{Z(z_0)<\infty}>0$.  By Lemma~\ref{howclose}, this implies that for
$\kappa\ge 8$ a.s.\ there is no $z\in\closure\H\setminus\{0\}$
which is swallowed.
\QED

\begin{exercise}Prove that $\Pb{\tau(z)<\infty}=0$
when $\kappa\in(0,4]$ and $z\in\H$.
\end{exercise}


\section{Transience}\label{trans}

\begin{theorem}[Transience]\label{transient}
For all $\kappa\ne8$ the \SLEk/ trace $\gamma(t)$ is
transient a.s.  That is,
$$\PB{\liminf_{t\to\infty} |\gamma(t)|=\infty}=1\,.$$
\end{theorem}

\begin{lemma}\label{avoid}
Suppose that $\kappa\le 4$, and let $x\in\R\setminus\{0\}$.
Then a.s.\ $x\notin\closure{\gamma[0,\infty)}$.
\end{lemma}
\proof
By symmetry, we may assume that $x>0$.  Moreover,
scaling shows that it is enough to prove the lemma in the
case $x=1$.  So assume that $x=1$.
Let $s\in(0,1/4)$ and suppose that there is some
first time $t_0>0$ such that $\bigl|\gamma(t_0)-1\bigr|=s$.
This implies that $\bigl|g_{t_0}(1/2)-\xo(t_0)\bigr|=O(s)$.
Indeed,
$\bigl|g_{t_0}(1/2)-\xo(t_0)\bigr|$
is the limit as $y\to\infty$ of $y$ times the harmonic
measure in $\H\setminus\gamma[0,t_0]$ from $i\,y$ of
the union of $[0,1/2]$ and the \lq\lq right hand side\rq\rq\
of $\gamma[0,t_0]$. By the maximum principle, this harmonic measure
is bounded above by the harmonic measure in $\H$ of the interval
$[1,\gamma(t_0)]$, which is $O(s/y).$

In the proof of Lemma~\ref{Rpart}, we have determined the probability
that $Y_x(t)=g_t(x)-\xo(t)$ hits $b$ before $a$, where $x\in[a,b]$ and $0<a<b$.
For $\kappa<4$, that probability goes to $1$ as $a\searrow 0$, uniformly
in $b$.  Consequently, a.s.\ there is a random $\eps>0$ such that for
all $b>x$, $Y_x(t)$ hits $b$ before $\eps$.  Since
$\sup_{t\in[0,T]} Y_x(t)<\infty$ a.s.\ for all finite $T$,
it follows that $\inf_{t\ge 0} Y_x(t)>0$ a.s.
Using this with $x=1/2$, we conclude from the previous paragraph
that a.s.\ $1\notin\closure{\gamma[0,\infty)}$, at least when $\kappa<4$.

Assume now that $\kappa=4$.  
Although the above proof does not apply (since $Y_x(t)$ is recurrent),
there are various different ways to handle this case.
Let $0<y<x<\infty$.  Let $\bar \gamma(t)$ denote the complex conjugate
of $\gamma(t)$, and set
$D_t:=\C\setminus\bigl(\gamma[0,t]\cup\bar\gamma[0,t]\cup (-\infty,y]\bigr)$.
Using Schwarz reflection, $g_t$ may be conformally extended to $D_t$.
Note that $g_t(D_t)=\C\setminus \bigl(-\infty,g_t(y)\bigr]$.
By the Koebe 1/4 Theorem applied to $g_t^{-1}$,
the distance from $x$ to $\p D_t$
is at least $\frac 14 \bigl(g_t(x)-g_t(y)\bigr)/g_t'(x)$.
Consequently, it suffices to show that a.s.
$$
\sup_{t\ge 0} g_t'(x)/\bigl(g_t(x)-g_t(y)\bigr)<\infty\,.
$$
Set $Q(t)=Q_{x,y}(t) := \log g_t'(x)-\log \bigl(g_t(x)-g_t(y)\bigr)$.
Then the above reduces to verifying that $\sup_{t} Q(t)<\infty$.
The time derivative of $Q$ is
\begin{equation*}
\p_t Q(t) =-2\,Y_x(t)^{-2} + 2 \,Y_x(t)^{-1}\,Y_y(t)^{-1}\,.
\end{equation*}
In particular, $Q(t)$ is increasing in $t$.
Define
$$
G(s):= \log s\, \log(1+s)-\frac 12 \log^2(1+s)+ \Li_2(-s)\,,
$$
where $\Li_2(z)=\int_z^0 s^{-1}\log(1-s)\, ds$.
Note that $G$ satisfies $s\,(1+s)^2\,G''(s)+s\,(1+s)\,G'(s)=1$.
A straightforward application of It\^o's formula shows that
$Q(t)-G_t$ is a local martingale, where
$$
G_t:=G\Bigl(\frac{g_t(x)-g_t(y)}{g_t(y)-\xo(t)}\Bigr)\,.
$$
Consequently, there is an increasing sequence of
stopping times $t_n$ tending to
$\infty$ a.s.\ such that
$\E[Q(t_n)] = \E[G_{t_n}]+Q(0)-G_0$.  It is immediate
to verify that $\sup\{ G(x):x>0\}<\infty$. 
Hence, $\limsup_{n\to\infty}\E[Q(t_n)]<\infty$.
Because $Q$ is increasing in $t$, it follows that
$\E[\sup_t Q(t)]<\infty$ and hence $\sup_t Q(t)<\infty$
a.s., which completes the proof.
\QED

\begin{lemma}\label{interior}
Suppose that $\kappa>4$, $\kappa\ne 8$, and let $t>0$.
Then a.s.~there is some $\epsilon>0$ such
that $K_t\supset\{z\in\H:|z|<\epsilon\}$.
\end{lemma}

\proof
Consider first the case $\kappa\in(4,8)$.
We know that $\gamma(t)$ is a continuous path, by Theorem~\ref{slecont}.
Note that by scale and reflection invariance, Lemma~\ref{ncard} shows that
for $a<b<0$ there is a positive probability that
the interval $[a,b]$ is swallowed all at once; i.e.,
$\Pb{[a,b]\cap\gamma[0,\infty)=\emptyset}>0$, but
a.s.\ there is some $t$ such that $[a,b]\subset K_t$.
When this happens, we have
$[a,b]\cap\closure{H_t}=\emptyset$.

Let $S$ be the set of limit points
of $g_{\tau(1)}(z)$ as $z\to0$ in $H_{\tau(1)}$. (If $0\notin \closure{H_{\tau(1)}}$, set
$S=\emptyset$.)
Then $S\subset \bigl(-\infty,\xo(\tau(1))\bigr)$.
Note that $\bigl( g_{\tau(1)}(\gamma(t+\tau(1)))-\xo(\tau(1)):t\ge 0\bigr)$
has the same distribution as $\bigl(\gamma(t):t\ge 0\bigr)$.
Consequently, by the above paragraph,
there is positive probability that
$g_{\tau(1)}(\gamma[\tau(1),\infty))$ is disjoint from $S$,
but there is some $t>\tau(1)$ such that
$g_{\tau(1)}(K_t\setminus K_{\tau(1)})\supset S$.
This implies that with positive probability there
is some $t>0$ such that
$$
\delta:=\Pb{0\notin\closure{H_t}}>0\,.
$$
By scale invariance, $\delta$ does not depend on $t$.
Since these events are monotone in $t$, it follows that
$$
\Pb{\forall t>0,\ 0\notin\closure{H_t}}=\delta\,.
$$
By Blumenthal's 0-1 law
it therefore follows that $\delta=1$.
This completes the proof in the case $\kappa\in(4,8)$.

Now suppose that $\kappa>8$.
The above argument implies that it is sufficient to show
that there is some $t>0$ such that with
positive probability $0\notin\closure{H_t}$.
Striving for a contradiction, we assume that for all
$t>0$ a.s.\ $0\in\closure{H_t}$.

Let $z_0\in\H$ be arbitrary.
We know from Theorem~\ref{tswal} and Lemma~\ref{posvol} that a.s.\
$z_0\in\gamma[0,\infty)$.  Set $t_0:=\tau(z_0)$.
Then $z_0=\gamma(t_0)$. Let $t>0$.
Conditioned on $\ev F_{t_0}$,
the set $g_{t_0}(K_{t+t_0}\setminus K_{t_0})$
has the same distribution as $K_t+\xo(t_0)$. It therefore follows
by our above assumption  that
a.s.\ $\xo(t_0)\in\closure {g_{t_0}(H_{t+t_0})}$.  Mapping by $g_{t_0}^{-1}$,
this implies that a.s.\ $z_0\in\closure {H_{t+t_0}}$; that is,
$z_0\in\partial K_{t+t_0}$.
We conclude that for every $z\in \H$ and every $t>0$ a.s.\
$z\in\partial K_{t+\tau(z)}$.
Therefore, a.s.\
$$
\area(\partial K_t)=\int_{z\in\H} 1_{z\in\partial K_t} \,dx\,dy
\ge \int_{z\in\H} 1_{t>\tau(z)}\,dx\,dy\,,
$$
and Fubini implies that $\partial K_t$ has positive measure with
positive probability for large $t$.
However, this contradicts Corollary~\ref{smallbd}, completing the proof
of the lemma.
\QED

\proofof{Theorem \ref{transient}}
First suppose that $\kappa>4$. Then we know from Lemma~\ref{interior}
that there is a random but positive $r>0$ such that
$K_1\supset \bigl\{z\in\H:|z|<r\bigr\}$.
Let $\epsilon>0$, and let $r_0=r_0(\kappa,\epsilon)>0$ be a constant such that
$\P[r<r_0]<\epsilon$.  Then the scaling property of SLE
shows that for all $t>0$
$$
\PB{K_t\supset \bigl\{z\in\H:|z|<r_0\sqrt t\bigr\}}\ge 1-\epsilon\,.
$$
This implies that
$$
\Pb{\exists t'>t\ : |\gamma(t')|<r_0\sqrt t}<\epsilon\,.
$$
Since $\epsilon>0$ was arbitrary, this clearly implies transience.

Consider now the case $\kappa\in[0,4]$.  In this range,
$\gamma$ is a simple path.  Then a.s.\ there are
two limit points $x_0,x_1$ for $g_1(z)$ as $z\to0$ in $H_1$.
Note that $g_1(\gamma[1,\infty))$ has the same
distribution as $\gamma[0,\infty)$ translated by $\xo(1)$.
Consequently, Lemma~\ref{avoid} shows that
a.s.\ $x_0$ and $x_1$ are not in the closure of
$g_1(\gamma[1,\infty))$.  This means that
$0\notin\closure{\gamma[1,\infty)}$ a.s.
The proof is then completed as in the case $\kappa>4$ above.
\QED

\begin{corollary}\label{spacefill}
Suppose that $\kappa>8$, then $\gamma[0,\infty)=\closure{\H}$
a.s.
\end{corollary}

\proof
{}From Lemma~\ref{howclose} and~\eref{dtz} we know
that $\gamma[0,\infty)$ is a.s.\ dense in $\closure{\H}$.
Since $\gamma$ is a.s.\ transient, $\gamma[0,\infty)$ is
closed in $\C$, and hence is equal to $\closure{\H}$.
\QED

\begin{update}
Corollary~\ref{spacefill} and Theorem~\ref{transient} are
true also for $\kappa=8$.
The proofs are based on the extension~\cite{\BLSWlesl} to $\kappa=8$
of Theorem~\ref{slecont}, and are otherwise the same.
\end{update}


\section{Dimensions}\label{dimensions}

In this section we study the size of the \SLE/ trace and the hull boundary.
When $\kappa\le 4$, these two sets are the same, and it is interesting
to compare the methods and results.
For sets $A\subset\C,$ denote by $N(\epsilon)=N(\epsilon,A)$ the minimal number of
disks of radius $\epsilon$ needed to cover $A$.

\subsection{The size of the trace}

\begin{theorem}\label{dimtrace} For every nonempty
bounded open set $A$ with $\closure{A}\subset\H,$
$$\limsup_{\epsilon\to 0} \frac {\log \EB{N\bigl( \epsilon, \gamma[0,\infty)
\cap A\bigr)}}{|\log \epsilon|} =
\begin{cases} 1+\kappa/8 & \kappa< 8,\\
2 & \kappa\geq 8\,.
\end{cases}
$$
\end{theorem}

\proof
We may assume $\kappa<8$ because $\gamma[0,\infty)=\overline \H$
if $\kappa>8$, by Corollary~\ref{spacefill},
and when $\kappa=8$ we know that for every
$z\in\H$ we have $z\in\gamma[0,\infty)$ a.s.,
by Theorem~\ref{tswal} and Lemma~\ref{posvol}.
Fix $\epsilon>0$ smaller than the distance
from $A$ to $\R$, and denote those points $z\in \epsilon \Z^2$ with
$\dist(z,A)<\epsilon$ by $z_j, 1\leq j\leq n$. Then $n$ is comparable to
$1/\epsilon^2$.
Set $N:=N\bigl( \epsilon, \gamma[0,\infty) \cap A\bigr)$. Obviously
$N\leq \#\{z_j : \dist(z_j,\gamma[0,\infty))<\epsilon\} \leq 100 N$.
{}From~\eref{dtz} and~\eref{comp} we get that
$\E[N]$ is comparable to
$$\sum_{j=1}^n \Pb{Z(z_j)>\epsilon^{-1}}.$$
Write $F(t)=\P[Z(i)>t]$.
Using~\eref{comp}, the proof of the theorem reduces to showing that
$$\limsup_{t\to\infty} \frac {\log F(t)}{\log t} = \kappa/8-1.$$
 By~\eref{gauss},
$$
\lim_{a\nearrow 1-\kappa/8} \;\frac{\hG_{a,\kappa}(1)}{1-\kappa/8-a} =
\frac{\sqrt\pi\, \Gamma(4/\kappa-1/2)}{2\, \Gamma(4/\kappa)}$$
and since~$\hG_{a,\kappa}(i)=1$ from Lemma~\ref{howclose} we obtain
$$C^{-1}\, \frac1{1-\kappa/8-a}\leq \Eb{Z(i)^a} \leq C\, \frac1{1-\kappa/8-a}$$
for all $a<1-\kappa/8$ sufficiently close to $1-\kappa/8$. Hence, for small $\Delta
=\Delta(\kappa)>0$,
$$
C\Delta\ge
\int_0^{\Delta}\Eb{Z(i)^{1-\kappa/8-\delta}}\,  \delta\, d\delta
=\int_0^{\Delta} \int_0^\infty (1-\kappa/8-\delta)\,
t^{-\kappa/8-\delta}\, F(t)\, dt\,  \delta\, d\delta
\,.
$$
{}From $\int_0^{\Delta} t^{-\kappa/8-\delta}\, \delta\, d\delta \ge C\, t^{-\kappa/8}/(\log t)^2$
for $t\ge e$ (where $\Delta$ is treated as a constant),
we obtain
$$\int_e^{\infty}\frac{t^{-\kappa/8}\,F(t)}{(\log t)^2}\, dt = O(1).$$
Using $F(t)/F(2t)=O(1)$ from~\eref{comp} this easily implies for large $t$
$$F(t)\le C\,t^{\kappa/8-1}\, (\log t)^2.$$
There is a sequence $t_n\to\infty$ with
$$F(t_n)\ge \frac{t_n^{\kappa/8-1}}{(\log t_n)^2},$$
for otherwise we would have
$\Eb{Z(i)^{1-\kappa/8}}=O(1)\, \int_e^\infty\frac{ dt}{t\,(\log t)^2}<\infty,$
contradicting  Lemma~\ref{howclose}.
\QED

An immediate consequence is
\begin{corollary}\label{Hautrace}
For $\kappa < 8,$ the Hausdorff dimension of $\gamma[0,\infty)$ is a.s.\ bounded above by
$1+\kappa/8$. \QED
\end{corollary}

\subsection{The size of the hull boundary}

When studying the size of the hull boundary $\p K_t$, it suffices
to restrict to $t=1$, by scaling.
We will first estimate the size of $\partial H_1\cap\{y>h\}$  for small $h>0$
in terms of the convergence exponent of the
Whitney decomposition of $H_1$. A Whitney decomposition of $H_1$ is a covering
by essentially disjoint closed squares $Q\subset H_1$ with sides parallel
to the coordinate
axes such that the side length $d(Q)$ is comparable to the distance of $Q$ from the
boundary of $H_1$. One can always arrange the side lengths
to be integer powers of 2 and that
\begin{equation}\label{whitney}
d(Q)\leq \dist(Q,\partial H_1)\leq 8\, d(Q).
\end{equation}
For example, one can take all squares of the form
$Q=2^n\bigl( [k,k+1]\times [j,j+1]\bigr)$
such that $Q\subset H_1$, $d(Q)\le \dist(Q,\p H_1)$, and
$Q$ is not strictly contained in any other $Q'$ satisfying these conditions.
One important feature of the Whitney decomposition is that the hyperbolic
diameters of the tiles in the decomposition are bounded away from $0$
and $\infty$.  (The hyperbolic metric on $H_1$ is the metric for
which the conformal map $g_1:H_1\to\H$ is an isometry, where
$\H$ is taken with the hyperbolic metric with length
element $y^{-1}\,|dz|$.)  This follows
from the well-known consequence of the Koebe 1/4 theorem
saying that the hyperbolic length element in
$H_1$ is comparable to $\dist(z,\p H_1)^{-1}\, |dz|$.
This description of the hyperbolic metric will be useful below.

Fix a Whitney decomposition of $H_1$, and denote by $W$ the collection of those
squares $Q$ for which $\dist(Q,K_1)\leq 1$ and $\sup\{\Im\,z:z\in Q\}\ge h$.
For $a>0$ let $$S(a)=S_h(a):=\sum_{Q\in W} d(Q)^a.$$

\begin{theorem}\label{whitneydim} 
Set
$$\delta(\kappa):=
\begin{cases} 1+\kappa/8 & \kappa\le 4,\\
1+2/\kappa & \kappa>4,
\end{cases}
$$
For $\kappa>0$, $\kappa\ne 4$, and all sufficiently small $h\in(0,1)$
we have
$$\E[S(a)]
\begin{cases} <\infty & a>\delta(\kappa),\\
= \infty & a\le\delta(\kappa). 
\end{cases}
$$
\end{theorem}

For $j,n\in\Z$, let
$x_{j,n}:=j\,2^{-n}$, $y_{n}:= 2^{-n}$ and
$z_{j,n} = x_{j,n}+i\,y_n$.
Let $I$ be the interval $g_1(\partial K_1)-\xo(1)$; that is, $\hf^{-1}(\p K_1)$.

\begin{lemma}\label{Wequiv} 
There are constants $C=C(h),c=c(h)>0$ such that the following holds.
Every $Q\in W$ is within distance $C$ in the hyperbolic
metric of $H_1$ from some $\hf_1(z_{j,n})$ with $x_{j,n}\in I$
and $n\ge 0$,
and every $\hf_1(z_{j,n})$ such that $x_{j,n}\in I$, $n\ge 0$ and
$\Im\, \hf_1(z_{j,n})\ge h$ is within hyperbolic distance $C$ from
some $Q\in W$.  Moreover, in each of these cases,
$$
c^{-1}\,d(Q)\le y_n\,\bigl|\hf_1'(z_{j,n})\bigr|\le c\,d(Q)\,.
$$
\end{lemma}
\proof
Observe that the diameter of $K_1$ is bounded from
below, because the uniformizing map
$g:\H\setminus K\to\H$ with the hydrodynamic
normalization~\eref{hydro} tends to
the identity as $K$ tends to a point,
and the uniformizing map $g_1$ for $K_1$
satisfies $g_1(z)=z+2/z+\cdots$ near $\infty$, by~\eref{power}.

Let $\iz\in Q\in W$.  It follows from the definition of $W$
and the fact that $\diam(K_1)$ is bounded from zero
that the harmonic measure of $K_1$ from $\iz$ in $H_1$
is bounded from zero.  Therefore,
the harmonic measure of $I$ from $\hf_1^{-1}(\iz)$
in $\H$ is bounded from zero.
This implies that there is some $z_{j,n}$ within
bounded distance from $\hf_1^{-1}(\iz)$ in the hyperbolic
metric of $\H$ such that $x_{j,n}\in I$.
Then $\hf_1(z_{j,n})$ is within bounded hyperbolic
distance from $Q$.
Because $\sup\{\Im\, z:z\in K_1\}=O(1)$,
we also have $\Im\, \iz=O(1)$,
and therefore $\Im\, \hf_1(z_{j,n})=O(1)$.
Since $\Im\, \hf_1(z_{j,n})\ge \Im\, z_{j,n}$,
it follows that $\max\{0,-n\}=O(1)$.
Then $z':=x_{j,n}+i\,2^{n\wedge 0}$
is equal to $z_{j',n'}$ for some $j'\in\Z,\, n'\in\N$,
and satisfies the requirements in the first statement of the lemma.

{}From the fact that $\diam(K_1)$ is bounded from zero
it also follows that $|I|$, the length of $I$, is bounded from zero.
Suppose that we have $x_{j,n}\in I$, $j\in\Z,\,n\in\N$.
Then the harmonic measure of $I$ from $z_{j,n}$
is bounded from zero.  Consequently,
the harmonic measure of $K_1$ from $\hf_1(z_{j,n})$ in
$H_1$ is bounded from zero.
Now suppose also that $\Im\,\hf_1(z_{j,n})\ge h$.
Since $\sup\{\Im\,z : z\in K_1\}=O(1)$ and
$\Im\,\hf_1(z_{j,n})=O(1)$, it follows
that the Euclidean distance from $\hf_1(z_{j,n})$ to $K_1$
is bounded from above.  If
$\dist\bigl( \hf_1(z_{j,n}),K_1\bigr)\le 1$, then
the Whitney tile $Q$ containing $\hf_1(z_{j,n})$
is in $W$.  Otherwise, let
$p$ be a point closest to $\hf_1(z_{j,n})$
among the points $z$ satisfying $\Im\,z\ge h$
and $\dist(z,K_1)\le 1$. Since $h<1$ and
$\dist\bigl(\hf_1(z_{j,n}),K_1\bigr)=O(1)$,
it follows that $\bigl|p-\hf_1(z_{j,n})\bigr|=O(1)$.
Since the points in $\bigl[p,\hf_1(z_{j,n})\bigr]$
have distance at least $h$ from $\p H_1$,
it follows that the hyperbolic distance from $p$
to $\hf_1(z_{j,n})$ is bounded above.  For the second statement
of the lemma we may then take $Q$ to be the
Whitney tile containing $p$.

The last statement holds because
$\hf_1$ is an isometry from the hyperbolic metric
of $\H$ to the hyperbolic metric of $H_1$.
\QED

\proofof{Theorem~\ref{whitneydim}}
Set
$$
\tilde S_h(a):=
\sum_{j=-\infty}^\infty\;\sum_{n=0}^\infty
1_ {x_{j,n}\in I}\,1_ {\Im \hf_1(z_{j,n})>h}\, |\hf_1'(z_{j,n})\, y_{n}|^a \,.
$$
We now show that this quantity is  comparable to $S(a)$.
Let $C$ be the constant in the lemma.
Since each square $Q\in W$ has a bounded diameter in the
hyperbolic metric of $H_1$, and the hyperbolic
distance in $\H$ between any two distinct $z_{j,n}$ is
bounded from below,
it follows that there are at most a bounded number of $\hf_1(z_{j,n})$
within distance $C$, in the hyperbolic metric of $H_1$, from any single $Q\in W$.
Similarly, there are at most a bounded number of $Q\in W$
within hyperbolic distance $C$ from any point in $H_1$.
Moreover, if $z$ is within hyperbolic distance
$C$ from some $Q\in W$,
then $\Im\,z$ is bounded away from zero,
because the hyperbolic distance in $H_1$ between any two points
is an upper bound on the distance between them in the hyperbolic
metric of $\H$ (by the Schwarz lemma, for instance).
Consequently, Lemma~\ref{Wequiv} implies
\begin{equation}\label{sabound}
c_1^{-1}\, \tilde S_h(a)
\le
S_h(a) \le c_1\, \tilde S_{h'}(a)
\,,
\end{equation}
for some constant $c_1,h'$,  which depends only on $h$ and $a$.
Now,
\begin{equation}\label{esa}
\Eb{\tilde S_h(a)}=
\sum_{n=0}^\infty \sum_{m=0}^\infty\sum_{j=-\infty}^{\infty}
\EB{ \bigl|\hf_1'(z_{j,n})\ y_{n}\bigr|^a \,
     1_{\ev A(j,n,m)} }
\end{equation}
where $\ev A(j,n,m)$ denotes the event $\{x_{j,n}\in I, \Im \hf_1(z_{j,n})\ge h, m\le |I|<m+1\}$.
We split this sum into two pieces, depending on whether $m\le n$ or $m>n$.
Clearly, $0\in I$.  It is easy to see that
$$
|I| \le O\Bigl(1+\max_{0\le t\le 1} \bigl|\xo(t)\bigr|\Bigr)\,,
$$
for example, by considering the differential equation~\eref{chordal} directly.
This implies $\Pb{|I|\ge m}\le O(1)\,e^{-C(\kappa)\,m^2}$,
for some $C(\kappa)>0$.
Since
$$
\bigl|y\,\hf'_1(x+i\,y)\bigr|\le O(1)\,\Im\, \hf_1(x+i\,y)\le O(1)
$$
for $0<y\le 1,$
we have the trivial estimate
\begin{align*}
\EB{\sum_{j=-\infty}^{\infty} |\hf_1'(z_{j,n})\ y_{n}|^a
     1_{\ev A(j,n,m)} }
&
 \leq O(1)\, m\, 2^{n}\, \Pb{|I|\ge m}
\\&
\le O(1)\, m\, 2^{n}\, e^{-C(\kappa) m^2}
\end{align*}
for all $m\ge0$.  This shows that in~\eref{esa} the terms
with $m\ge n$ have finite sum.
For $m\le n$, we use Theorem~\ref{derexp} together with~\eref{nearhight} and obtain
$$
\EB{\sum_{j=-\infty}^{\infty} \bigl|\hf_1'(z_{j,n})\ y_{n}\bigr|^a \,
     1_{\ev A({j,n,m})} }
 \leq C(h,\kappa,b)\, m^{2b+1}\, 2^{-n\,(\lambda-2b-1)},
$$
where $a,\lambda,b$ and $\kappa$ are related by~\eref{al} with $\nu=1$ and $b>0$.
For $b$ between $1/2$ and $2/\kappa$ we have $\lambda-2b-1>0$ and conclude
that $ \Eb{\tilde S(a)}<\infty $ in this range. By~\eref{sabound}, $\E[S(a)]<\infty$
for $\kappa< 4, a> 1+\kappa/8$ and for $\kappa> 4, a>1+2/\kappa$.

Now consider $a=\delta(\kappa)$.
Suppose that $2^{-n}<h/4$ and $\hf_1(z_{j,n})\ge h$. 
We claim that the harmonic measure in $H_1$ of $K_1$ from $\hf_1(z_{j,n})$
is bounded from zero.  Indeed, 
fix a constant $h_1>1$ such that
$h_1/2\ge \sup\{\Im\,z:z\in K_1\}$.  (It is easy to see that $h_1=4$ works.)
Consider the harmonic function $Y(z):=\Im\, g_1(z)$ on $H_1$.
Then $Y(x+i\,y)-y\to 0$ uniformly in $x$, as $y\to\infty$, by
the hydrodynamic normalization~\eref{hydro} of $g_1$.
Since $\H+\frac i2 \,h_1\subset H_1$,
the maximum principle applied to $Y(z)-\Im z$ implies that $Y(x+i\,h_1)\ge h_1/2$ for
all $x\in\R$.  Consequently, for every $z\in H_1$ with
$\Im\,z< h_1$, the harmonic measure in $H_1\cap\{y<h_1\}$ from $z$ of
the line $y=h_1$ is at most $2\,Y(z)/h_1$.
However, the harmonic measure from $z$
of $y=h_1$ in the domain $\H\cap \{y<h_1\}$ is $\Im\,z/h_1$.
Therefore, the harmonic measure from $z$ of $K_1$
is at least the difference between these two quantities,
which is at least $\Im\,z/h_1- 2\, Y(z)/h_1$.
Applying this to $z:=\hf_1(z_{j,n})$ proves our claim
that the harmonic measure of $K_1$ from $\hf_1(z_{j,n})$ in $H_1$
is bounded from zero.  This implies that the
harmonic measure of $I$ from $z_{j,n}$ in $\H$ is bounded from
zero, which implies that at bounded hyperbolic distance from
$z_{j,n}$ we can find a point $z_{j',n'}$ with
$x_{j',n'}\in I$.
Consequently, to prove that $\Eb{\tilde S_h\bigl(\delta(\kappa)\bigr)}=\infty$
for all sufficiently small $h>0$, it suffices to prove
that for all $n_0\in\N$ and all sufficiently small $h>0$
$$
\EB{
\sum_{j=-\infty}^\infty\;\sum_{n=n_0}^\infty
1_ {\Im \hf_1(z_{j,n})>h}\, |\hf_1'(z_{j,n})\, y_{n}|^{\delta(\kappa)}}=\infty \,.
$$
This follows from Lemma~\ref{lbdexp} with $b=1/2$ if $\kappa<4$
or $b=2/\kappa$ if $\kappa>4$.  An appeal to~\eref{sabound}
gives $\Eb{S\bigl(\delta(\kappa)\bigr)}=\infty$.
Since $\Eb{S(a)}<\infty$ for $a>\delta(\kappa)$,
it also follows that $\Eb{S(a)}=\infty$ for $a<\delta(\kappa)$,
and the proof is complete.
\QED

\begin{corollary}\label{lowerboxdim} 
For all $\kappa>0$, $\kappa\ne4$, and all sufficiently small $h>0$
$$\limsup_{\epsilon\to0}\;
 \bigl|\log\eps\bigr|^{1+h}\,\epsilon^{\delta(\kappa)}\,
\EB{N\bigl(\epsilon, \partial K_1\cap \{y>h\}\bigr)}= \infty.
$$
\end{corollary}
\proof
Denote by $W(n)$ the number of $Q\in W$ with $2^{-n-1}< d(Q)\le 2^{-n}$. Then~\eref{whitney}
immediately gives
\begin{equation}\label{wversusn}
W(n)\leq C\, N(2^{-n}, \partial K_1\cap \{y>h\}).
\end{equation}
The corollary now follows from
$
S(a)=\sum_{n=0}^\infty W(n) 2^{-a n}
$
and Theorem~\ref{whitneydim}.
\QED

Finally we state an upper bound for the covering dimension.

\begin{theorem}\label{upperboxdim}
For all $\kappa,h>0$  and $a>\delta(\kappa),$ a.s.\
$$\lim_{\epsilon\to0}
\eps^a\,N\bigl(\epsilon, \partial K_1\cap \{y>h\}\bigr)= 0.
$$
\end{theorem}

\proof
Suppose first that $\kappa\ne 4$.
By Theorem~\ref{holdercont} we know that $\hf_1$ is a.s.\ H\"older continuous.
There is no lower bound version of~\eref{wversusn} for general domains.
However, by results of Pommerenke~\cite{\BPomm} and
Makarov~\cite{\BMakFine} it is known that for H\"older domains
(i.e., the conformal map from $\U$ into the domain is H\"older)
the asymptotics of $N$ can be recovered from the
convergence exponent of the Whitney decomposition.
Our proof follows this philosophy.
We claim that there is a (random) constant $C=C(\omega)>0$, depending on $h$
and the H\"older norm of $\hf_1$, such that for every
small $\eps>0$ and every point $z_0\in\p K_1\cap\{y\ge 2 h\}$ 
\begin{equation}\label{e.KR} 
\exists Q\in W,\qquad 
Q\subset D(z_0,\eps),\;\;\eps\ge d(Q)\ge C \eps/|\log\eps|\,.
\end{equation}
In fact, this follows from~\cite{\BKR}, but we include a short proof 
here, for the convenience of the reader.

Let $x_0\in\R$ be a point such that $\hf_1(x_0)=z_0$.
Clearly, there is such an $x_0$, because $\hf_1$ extends continuously to $\closure{\H}$. 
Let $\beta(s):= \hf_1(x_0+i\,s)$, $s>0$,
let $\hat Q(s)$ be the set of Whitney cells
$Q\in W$ which intersect $\beta\bigl((0,s)\bigr)$, and
set 
$$
q(s):= \sum_{Q\in\hat Q(s)} d(Q)\,.
$$
Since $\hf_1$ is an isometry between
the hyperbolic metric of $\H$ and that of $H_1$,
it follows that for every $s>0$ the arc
$\beta[s,2s]$ meets a bounded number of Whitney cells,
and each of these has Euclidean diameter comparable to
the Euclidean diameter of $\beta[s,2s]$.  Suppose
that $\hf_1$ is H\"older with exponent $\alpha>0$.
Then, $d(Q)\le O(1)\,s^{\alpha}$ for $Q\in\hat Q(s)$,
and consequently also $q(s)\le O(1)\, s^\alpha$.
Observe that there is a constant $c_0$ so that
$q(s)/\sup\bigl\{q(s'):s'<s\bigr\}<c_0$,
because there are constant upper bounds on the number
of Whitney cells containing a point, and on
the ratio $d(Q)/d(Q')$ for neighboring
Whitney cells $Q$ and $Q'$. 
Let $s_0:=\sup\{s:q(s)<\eps\}$ and $s_1:=\sup\{s:q(s)<c_0^{-1}\,\eps/2\}$.
The choice of $c_0$ ensures that $q(s_0)-q(s_1)\ge \eps /2$.
Clearly $\hat Q(s_0)\setminus \hat Q(s_1)$
has at most $O\bigl(1+ \log(s_0/s_1)\bigr)$ cells.
Therefore, there is some $Q\in\hat Q(s_0)\setminus \hat Q(s_1)$
with $d(Q)\ge \eps/O\bigl(1+ \log(s_0/s_1)\bigr)$.
Since $q(s)\le O(1)\,s^\alpha$, it follows that $O(1)\,s_1\ge \eps^{1/\alpha}$,
and we get $O(1)\,d(Q)\ge \eps/|\log\eps|$, proving~\eref{e.KR}.

We know that for every set $F\subset\C$ the covering number
$N(\eps,F)$ is bounded by a constant times
the maximum number of disjoint disks of radius $\eps$ centered in $F$.
Therefore, by~\eref{e.KR}, $N\bigl(\eps, \partial K_1\cap \{y>2\,h\}\bigr)$
is bounded by a constant times the number of Whitney cells in $W$
with size between $O(1)\,\eps$ and $\eps/O\bigl(|\log\eps|\bigr)$; that is,
$$
N\bigl(2^{-n}, \partial K_1\cap \{y>2\,h\}\bigr)
\le C(\omega)\,\sum_{j=n}^{n+O(\log n)} W(j)\,,
$$
which implies for all $a'>0$,
$$
N\bigl(2^{-n}, \partial K_1\cap \{y>2\,h\}\bigr)
\le C(\omega)\, 2^{(n+O(\log n))\,a'} S(a')\,.
$$
By Theorem~\ref{whitneydim}, $S(a')<\infty$ a.s.\ for $a'>\delta(\kappa)$.
Therefore, the case $\kappa\ne 4$ in Theorem~\ref{upperboxdim} now follows
by taking $a'\in\bigl(\delta(\kappa),a\bigr)$.


The following argument covers not only the case $\kappa=4$, but the
whole range $\kappa\in(0,4]$.
For such $\kappa$, Theorem~\ref{dimtrace} can be used,
since $\p K_1=\gamma[0,1]$ for $\kappa\in(0,4]$.
Although Theorem~\ref{dimtrace} discusses the expected
number of balls needed to cover $\gamma[0,\infty)\cap A$,
where $A$ is an arbitrary bounded set satisfying $\closure{A}\subset\H$,
one can refine the result to bound $N\bigl(\epsilon,\p K_1\cap\{y>H\}\bigr)$,
as follows.  Given a point $z=x+i\,y$, $y>h$, we need to estimate
the probability that $z$ is close to $K_1$.  We may certainly
restrict attention to the case where $y$ is bounded by some
constant.  There is a constant $C\ge 1$ such that the probability that $z$
is close to $K_1$ is
bounded by the probability that there is some time $t_1<1$ such that
$\xo(t_1)$ is within distance $C$ from $x$ times the
probability that $Z(z)$ is large, conditioned on the existence
of such a $t_1$.  The proof is then completed by applying
the estimates in the proof of Theorem~\ref{dimtrace} and
the Markov property.  The details are left to the reader.
\QED

The upper bound for the Hausdorff dimension follows immediately:

\begin{corollary}\label{Haubd}
For $\kappa > 4,$ the Hausdorff dimension of $\partial K_1$ is
a.s.\ at most $1+2/\kappa$. \QED
\end{corollary}

\section{Further discussion and some open problems}\label{further}

We believe Theorem~\ref{slecont} to hold also in the case $\kappa=8$:

\begin{conjecture}\label{slecont8}
The \SLEkk8/ trace is a continuous path a.s.
\end{conjecture}

\begin{update}
This is proved in~\cite{\BLSWlesl}.
\end{update}

\begin{conjecture} The Hausdorff dimension of the \SLEk/ trace is
a.s.\ $1+\kappa/8$ when $\kappa<8$.
\end{conjecture}

\begin{update}
A proof of this result is announced in~\cite{\BBeffaraHaus}
for $\kappa\ne 4$.
\end{update}

Corollary~\ref{Hautrace} shows that the Hausdorff dimension cannot
be larger than $1+\kappa/8$.

\begin{conjecture}[Richard Kenyon] The Hausdorff dimension of $\p K_1$ is
a.s.~$1+2/\kappa$ when $\kappa\ge 4$.
\end{conjecture}

Corollary~\ref{Haubd} shows that the Hausdorff dimension cannot be
larger than $1+2/\kappa$.

\begin{problem} Find an estimate for the modulus of continuity for
the trace of \SLEkk4/.  Improve the estimate for the modulus
of continuity for \SLEk/ for other $\kappa$.
What better modulus of continuity can be obtained with a different
parameterization of the trace?
\end{problem}

Based on the conjectures for the dimensions~\cite{\BDup},
Duplantier suggested that \SLEkk{16/\kappa}/ describes the boundary
of the hull of \SLEk/ when $\kappa>4$.  This has further support,
some of which is discussed below.
However, at this point there isn't even a precise version for
this duality statement.

\begin{problem}\label{Dup} Understand the relation between the boundary
of the \SLEk/ hull and \SLEkk{16/\kappa}/ when $\kappa>4$.
\end{problem}

\begin{update}
See~\cite{\DubedatDuality} for recent progress on this problem.
\end{update}

We now list some conjectures about discrete processes tending to \SLEk/
for
various $\kappa$.
Let $D$ be a bounded simply connected domain in $\C$ with smooth
boundary,
and let $a$ and $b$ be two distinct points in $\p D$.  Let $\phi$ be
a conformal map from $\H$ to $D$ such that $\phi(0)=a$ and
$\phi(\infty)=b$.
The {\bf \SLE/ trace in $D$ from $a$ to $b$} is defined as the image
under $\phi$ of the trace of chordal \SLE/ in $\H$.
Similarly, one defines radial \SLE/ in arbitrary simply
connected domains strictly contained in $\C$.

\begin{figure}
\centerline{\includegraphics*[height=2.3in]{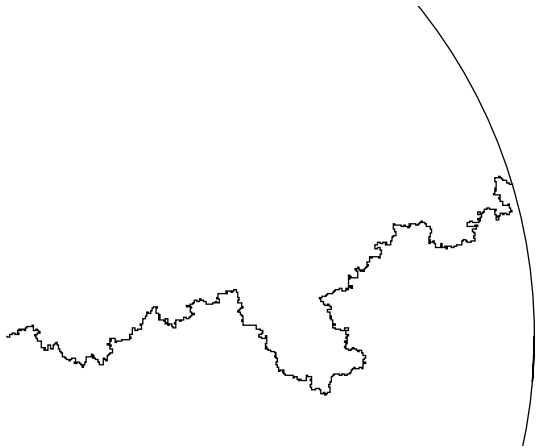}}
\caption{\label{f.lerw}A sample of the loop-erased
random walk; believed to converge to radial \SLEkk2/.}
\end{figure}

In~\cite{\BSchSLE} it was shown that if the loop-erased
random walk in simply connected domains in $\C$ has a conformally
invariant scaling limit, then this must be radial \SLEkk2/.
See Figure~\ref{f.lerw} for a sample of the loop-erased random
walk in $\U$.

\begin{figure}
\centerline{\includegraphics*[height=2.1in]{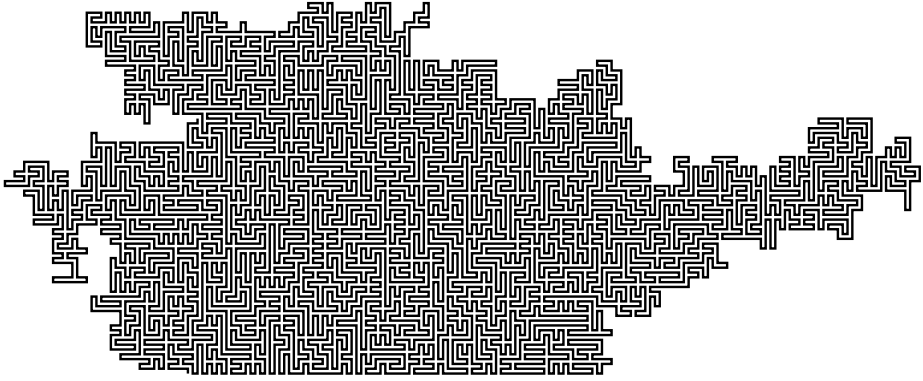}}
\caption{\label{f.peano}A piece of the UST Peano curve;
believed to converge to \SLEkk8/.}
\end{figure}

In~\cite{\BSchSLE} it was also conjectured that the UST Peano curve,
with appropriate boundary conditions (see Figure~\ref{f.peano}),
tends to \SLEkk8/ as the mesh goes to zero.
This offers additional
support to Duplantier's conjectured duality~\ref{Dup},
Since the outer boundary of the hull of the UST Peano curve
consists of two loop-erased random walks, one in the tree and
one in the dual tree.  See~\cite{\BSchSLE} for additional details.

\begin{update}
The convergence of LERW to
\SLEkk2/ and of the UST Peano path to \SLEkk8/ is proved in~\cite{\BLSWlesl}.
\end{update}

\begin{conjecture} Consider a square grid of mesh $\epsilon$
in the unit disk $\U$.  Then the uniform measure on simple
grid paths from $0$ to $\p\U$ on this grid converges weakly
to the trace of radial \SLEkk8/ as $\eps\searrow0$, up
to reparameterization.
The uniform measure on simple paths
joining two distinct fixed boundary points converges weakly to
the trace of chordal \SLEkk8/, up
to reparameterization.
\end{conjecture}

We have not bothered to make this conjecture entirely precise,
since there is more than one natural choice for the topology
on path space.  However, this does not seem to be very
important at this point.

Note that we do not require that the grid path be space 
filling (i.e., visit every vertex).  Indeed, it is not hard to
see that a uniformly selected path will be essentially filling,
in the sense that the probability that a nonempty open set inside
$\U$ is unvisited goes to zero as $\eps\searrow0$.

\medskip

\begin{figure}
\centerline{\includegraphics*[width=2.0in]{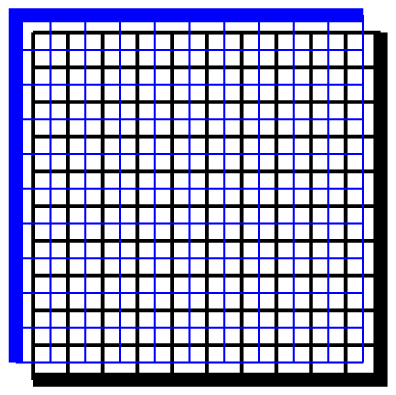}}
\caption{\label{f.dua}A partially wired grid and its dual.}
\end{figure}

Consider an $n\times n$ grid in the
unit square, and identify the vertices of the right
and bottom boundary arcs to a single vertex $w_0$.
Let $G_n$ denote this graph, and note that the
dual graph $G_n^\dagger$ is obtained by
rotating $G_n$ 180$^o$ and shifting appropriately.
See Figure~\ref{f.dua}.
Let $\omega$ be a sample from the
critical random-cluster
measure on $G_n$ with parameters $q$ and $p=p(q)=\sqrt q/(1+\sqrt q)$
(which is the conjectured value of $p_c=p_c(q)$)
and let $\omega^\dagger$ be the set of edges in $G_n^\dagger$
that do not cross edges of $\omega$.
Then the law of $\omega^\dagger$ is the random-cluster measure on
$G_n^\dagger$ with the same parameters as for $\omega$.
Let $\beta$ be the set of points in the unit square
that are at equal distance from the component of $\omega$ containing
$w_0$ and the component of $\omega^\dagger$ containing the
wired vertex of $G_n^\dagger$.

\begin{conjecture}
If $q\in(0,4)$, then
as $n\to\infty$, the law of $\beta$ converges weakly
to the trace of \SLEk/, up to reparameterization,
where
$ \kappa=\frac{4\pi}{\cos^{-1}\bigl(-\sqrt{q}/2\bigr)} $.
\end{conjecture}

The above value of $\kappa$ is derived by applying the
predictions \cite{\BDupSal} of the dimensions of
the perimeters of clusters in the random cluster model.

\begin{figure}
\centerline{\includegraphics*[width=4.0in]{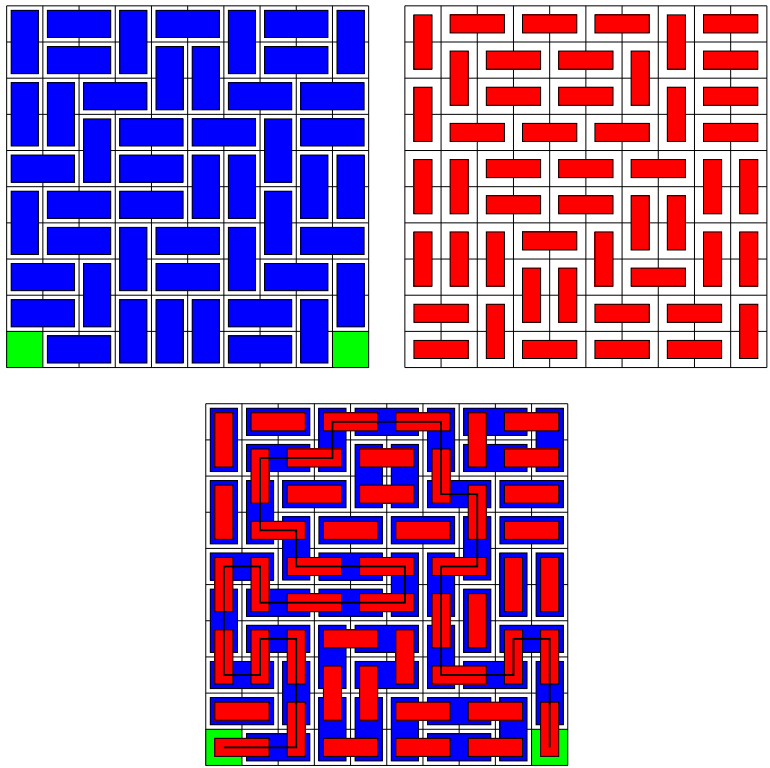}}
\caption{\label{f.ddom}The double domino path.}
\end{figure}

Consider a random-uniform domino tiling of the unit square, such as
in the top right of Figure~\ref{f.ddom}.
Taking the union with another independent
domino tiling with the two lower corner squares removed,
one gets a path joining the two removed squares; see
Fig.~\ref{f.ddom}.  This model was
considered by~\cite{\BRHAdd}.  Richard Kenyon proposed the following

\begin{problem} Does this double-domino model converge weakly
(up to reparameterization) to \SLEkk4/
as the mesh of the grid tends to zero?
\end{problem}

Kenyon also produced some calculations for this model that
agree with the corresponding calculations for \SLEkk4/.

\medskip

Let $a\in\R$, $\eps>0$, and
consider the following random walk $S$ on the grid
$\epsilon\,\Z^2$.
Set $S_0=0$, and inductively, we construct $S_{n+1}$ given
$\{S_0,\dots,S_n\}$. If $S_n\notin\U$, then stop.  Otherwise,
let $h_n:\eps\,\Z^2\to[0,1]$ be the function which is $1$ outside of $\U$,
$0$ on $\{S_0,\dots,S_n\}$, and discrete-harmonic elsewhere.
Now choose $S_{n+1}$ among the neighbors $v$ of $S_n$ such that $h(v)>0$,
with probability proportional to $h(v)^a$.
When $a=1$, this is called {\bf Laplacian random
walk}, and is known to be equivalent to loop-erased random walk.
When $a=0$ this walk is believed to be closely related to
the boundary of percolation clusters.

The following problem came up in discussions with Richard Kenyon.

\begin{problem}  For what values of $a$ does this process converge
weakly to radial \SLE/?  What is the correspondence $\kappa=\kappa(a)$?
\end{problem}

Motivated in part by Duplantier's duality conjecture~\ref{Dup},
in~\cite{\BLSWrest} it is conjectured that chordal \SLEkk{8/3}/
is the scaling limit of the self-avoiding walk in a half-plane.
(That is, the limit as $\delta\searrow 0$ of the
limit as $n\nearrow\infty$ of the law of a random-uniform $n$-step simple
path from $0$ in $\delta\Z^2\cap\closure{\H}$.)
Some strong support to this conjecture is provided there,
which in turn can be viewed as further support for the duality in~\ref{Dup}.

\medskip
Many of the processes believed to converge to chordal \SLE/ have
the following reversibility property.  Up to reparameterization,
the path from $a$ to $b$ in $D$ ($a,b\in\p D$) has the same
distribution as the path from $b$ to $a$.  
It is then reasonable to conjecture:

\begin{conjecture} The chordal
\SLE/ trace is reversible for $\kappa\in[0,8]$.
\end{conjecture}

It is not too hard to see that this fails for $\kappa>8$.  In particular,
when $\kappa>8$ the chordal \SLE/ path is more likely to visit $1$ before
visiting $i$.  Since inversion in the unit circle preserves $1$ and $i$,
the probability for visiting $i$ before $1$ would have to be $1/2$ if one
assumes reversibility.  To apply this argument one has to note that a.s.\
the path visits $1$ and $i$ exactly once.

\begin{update}
The above conjecture is known to hold when $\kappa=0,2,8/3,6,8$.
The cases $\kappa=2,8$ follow from~\cite{\BLSWlesl} (though some
additional work is needed in the case where $\kappa=2$, since
the convergence of the LERW is to radial \SLE/), the case $\kappa=6$
follows from Smirnov's Theorem~\cite{\BSmirnovPerc},
and the case $\kappa=8/3$ follows from the characterization
in~\cite{\BLSWrest} of
\SLEkk{8/3}/ as the restriction measure with exponent $5/8$.
\end{update}

\medskip

Some conjectures and
simulations concerning the winding numbers of
\SLE/'s are described in~\cite{\BWeilWil}.

{\bigskip\noindent\bf Acknowledgments}.
The authors thank Richard Kenyon for useful discussions and Jeff Steiff for
numerous comments on a previous version of this paper.

\bibliographystyle{halpha}
\bibliography{mr,prep,notmr}

\end{document}